\documentclass[a4paper,12pt]{article}

\usepackage{authblk}
\usepackage[noadjust]{cite}
\usepackage{amsmath}
\usepackage{amsthm,amsfonts,amssymb}           
\usepackage{bbm}
\usepackage{mathtools}\mathtoolsset{showonlyrefs=true}
\usepackage[left=2.5cm,
            right=2.5cm,
            top=2.5cm,
            bottom=4cm
	       ]{geometry}
\usepackage[pdfstartview=FitH,
            CJKbookmarks=true,
            bookmarksnumbered=true,
            bookmarksopen=true,
            bookmarksdepth=2,
            colorlinks,
            linkcolor=blue,
            anchorcolor=blue,
            citecolor=blue,
            urlcolor=blue
            ]{hyperref}
\newtheorem{mthm}{Theorem}
\newtheorem{thm}{Theorem}[section]

\newtheorem{lem}[thm]{Lemma}

\theoremstyle{remark}

\def\R{{\mathbb R}}

\def\bd{\operatorname{bd}}
\def\relint{\operatorname{relint}}
\def\intr{\operatorname{int}}
\newcommand{\GL}{\mathrm{GL}}
\newcommand{\SL}[1]{\mathrm{SL}({#1})}                  
\newcommand\slpm[1]{\mathrm{SL}^\pm({#1})}
\newcommand{\GLn}{\mathrm{GL}(n)}                  
\newcommand\gln{\operatorname{GL}(n)}
\newcommand{\SLn}{\mathrm{SL}(n)}                  
\newcommand{\sln}{\mathrm{SL}(n)}                  
\newcommand{\Sn}{S^{n-1}}                  

\def\MP{\mathcal{P}}

\def\MN{\mathcal{N}}

\def\MH{\mathcal{H}}

\def\MPoon{\MP_{(o)}^n}

\def\MPon{\MP_{o}^n}
\def\MPo2{\MP_{o}^2}
\newcommand{\mpn}{\MP^n}

\def\sgn{\operatorname{sgn}}    
\def\a{\alpha}

\newcommand\ab[1]{\left(#1\right)}        
\newcommand\abs[1]{\left|#1\right|}        
\newcommand\abb[1]{\left\{#1\right\}}        

\newcommand{\cross}[1]{#1^{\times}} 
\newcommand{\ct}[3]{#1_{#2}^\times({#3})} 
\newcommand{\tenset}[2]{(\R^{#1})^{\otimes_{#2}}}
\newcommand{\ten}[1]{\mathbbm{t}}  
\newcommand{\symtp}[2]{{#1}^{#2}}
\newcommand{\symtenset}[2]{\operatorname{Sym}^{#2}(\R^{#1})}

\newcommand{\sympro}{\odot}
\newcommand{\symg}[1]{\operatorname{S}_{#1}}
\newcommand{\alttenset}[2]{\operatorname{Alt}^{#2}(\R^{#1})}
\newcommand{\dten}[1]{A^{(#1)}}
\newcommand{\tenval}{\mathrm{TVal}}

\newcommand*{\dif}{\mathop{}\!\mathrm{d}}             

\newcommand{\polyeo}{\mathcal{P}^2_o}	

\newcommand{\seq}[2]{\{#1,\dots,#2\}}  
\newcommand{\seqnb}[2]{#1,\dots,#2}  
\newcommand{\seqalpha}[3]{#1_{#2},\dots,#1_{#3}}  
\newcommand{\relseq}[4]{#1_{#2} #4 \dots #4 #1_{#3}} 
\newcommand{\set}[1]{\left\{#1\right\}}

\newcommand{\lct}{\,\mathcal{E}} 
\newcommand{\indval}[3]{Z_{#1;#2;#3}}
\newcommand{\indf}[3]{#1;#2;#3}
\newcommand{\indfconst}[2]{#1;#2}

\newcommand{\inpd}[2]{\langle #1,#2 \rangle} 
\newcommand{\polyn}{\MP^n}
\newcommand{\polye}{\MP^2}

\newcommand{\spl}{\mathcal T^n_o}					
\newcommand{\splo}{\mathcal T^n_o\setminus\set{o}}	

\makeatletter
\newcommand{\subjclass}[2][1991]{%
  \let\@oldtitle\@title%
  \gdef\@title{\@oldtitle\footnotetext{#1 \emph{Mathematics subject classification.} #2}}%
}
\newcommand{\keywords}[1]{%
  \let\@@oldtitle\@title%
  \gdef\@title{\@@oldtitle\footnotetext{\emph{Key words and phrases.} #1.}}%
}
\makeatother

\title{\bf{SL($n$) Contravariant Tensor Valuations of Orders up to $n$}}
\author[1]{Jin Li}
\author[2]{Dan Ma}
\affil[1]{Department of Mathematics and Newtouch Center for Mathematics, Shanghai University, Shanghai, China, 200444 \href{mailto:li.jin.math@outlook.com}{li.jin.math@outlook.com}}

\affil[2]{Department of Mathematics, Shanghai Normal University, Shanghai 200234, China, \href{mailto:madan@shnu.edu.cn}{madan@shnu.edu.cn}}
\date{}

\begin{document}

\maketitle

\begin{abstract}
A complete classification of \(\mathrm{SL}(n)\) contravariant, \(p\)-order tensor valuations on polytopes in \( \mathbb{R}^n \) for \( n \geq p \) is established without imposing additional assumptions, particularly omitting any symmetry requirements on the tensors. Beyond recovering known symmetric tensor valuations, our classification reveals asymmetric counterparts associated with the cross tensor and the Levi-Civita tensor. Additionally, some Minkowski-type relations for these asymmetric tensor valuations are obtained, extending the classical Minkowski relations of surface area measures.
\end{abstract}

\section{Introduction}
Let $\mathcal{Q}^n$ be a collection of subsets in $\R^n$.
A function $Z$ defined on $\mathcal{Q}^n$ and taking values in an abelian semigroup is called a valuation if
\begin{equation}\label{eqn:val0}
	Z(P)+Z(Q)=Z(P\cup Q)+Z(P\cap Q)
\end{equation}
for every $P,Q \in\mathcal{Q}^n$ such that $P\cap Q, P\cup Q \in \mathcal{Q}^n$.
We say a valuation is \emph{simple} if it vanishes on lower-dimensional sets.
The concept of valuation originated from Dehn's solution of Hilbert's Third Problem. Specifically, Dehn's solution demonstrated that volume is not the only real-valued, simple valuation invariant under rigid motions, thereby revealing the non-triviality of classifying rigid motion invariant valuations on polytopes (convex hull of finitely many points), a problem that remains open to this day. Subsequent developments, guided by Felix Klein's Erlangen program, have advanced the classification of valuations compatible with transformations such as rigid motions and affine transformations. A landmark achievement in this direction is Hadwiger's classification theorem \cite{Had57}, which fully characterizes continuous, rigid motion invariant, real-valued valuations on convex bodies as linear combinations of intrinsic volumes. This framework has since served as a cornerstone in integral geometry and geometric measure theory; see the books and surveys \cite{KR97,JK17,MR4703948,MR3380549}.
Recently, Hadwiger's theorem on convex functions was established by Colesanti, Ludwig, and Mussnig in \cite{CLM2020hadwiger1,CLM2025hadwiger2,CLM2022hadwiger3,CLM2023hadwiger4}.
Connections between valuation theory and differential geometry have been explored widely; see, for example, \cite{AB2017con,MR2746948}.

A tensor valuation is a valuation taking values in the space of tensors. In this paper, we study tensor valuations compatible with the special linear group of $\R^n$.
This places the problem within affine valuation theory and convex geometry.
Studying geometric functionals compatible with affine transformations is one of the central topics in convex geometry; see, e.g. \cite{LYZ00a,MR4618955,HP14b,LR10,MiP,BL2017Min,GHW13,Lut96}.

Let $\tenset{n}{p}$ be the set of all $p$-order tensors on $\R^n$ for $p \ge 0$,
and let $\symtenset{n}{p}\subset\tenset{n}{p}$ be the subset of symmetric ones.
We understand $\tenset{n}{0} = \symtenset{n}{0}=\R$.
First examples of affine tensor valuations are $M^{p,0},M^{0,p}: \MPoon\to\symtenset{n}{p}$ ($n\geq2, ~p\geq0$) defined by
\begin{equation}\label{eqn:tenp0}
	M^{p,0}(P):=\int_Px^p \dif x,  ~P \in \MPoon
\end{equation}
and
\begin{equation}\label{eqn:ten0p}
	M^{0,p}(P):=\int_{\Sn} h_P^{-p}(u) \symtp{u}{p} \dif V_P(u), ~P \in \MPoon.
\end{equation}
Here $\symtp{x}{p}$ is the $p$-fold tensor power of $x \in \R^n$,
$\MPoon$ is the space of polytopes in $\R^n$ that contain the origin in their interiors,
$\Sn\subset\R^n$ is the unit sphere, $V_P$ is the cone-volume measure of $P$ and $h_P$ is the support function of $P$; see definitions in \S \ref{sec:preliminaries}.
When $p=0$, both operators agree with the volume.
Within the framework of centro-affine Hadwiger theory, $M^{p,0}$ and $M^{0,p}$ are characterized through a series of seminal contributions by Haberl and Parapatits \cite{HP14b,HP16mom,HP17TV}, building on the foundational framework established in Ludwig's pioneering investigations \cite{Lud02a,Lud02c,Lud03} and Ludwig and Reitzner \cite{LR10}.
See \cite{Lud11b,Lud13} for some extensions to valuations on functions.

The tensors $M^{p,0},M^{0,p}$ are closely related to Brunn-Minkowski Theory and its $L_p$ extensions.
Minkowski problems for the cone-volume measures and $L_p$ surface area measures are one of the major problems in convex geometry and geometric PDE; see, for example, \cite{Lut93,BLYZ13,HLYZ2016dual,MR423267,MR4764744,MR2254308,BHO25}.
The operator $M^{p,0}$ belongs to the class of Minkowski tensors compatible with affine transformations (specifically, it is $\sln$ covariant), while $M^{0,p}$ belongs to the 
$L_p$ version version of that class of tensors (and is $\sln$ contravariant).
See \cite{McM97,Ale99b,Ale99,BDS21,BH18,HS14} for characterizations of Minkowski tensors compatible with rigid motions and \cite{ABDK19,LS17,FLR24} for other related tensor valuations.
For $p=2$, $M^{2,0}(P)$ and $M^{0,2}(P)$ are positive definite matrices and their associated ellipsoids are the Legendre ellipsoid and the Lutwak-Yang-Zhang ellipsoid of $P$, respectively; see, for example, \cite{MiP,LYZ00b,Lud03}.
For general integer $p>0$, we view $M^{p,0}$ and $M^{0,p}$ as $p$-linear, symmetric forms on $\R^n$ given by
\begin{align}
	M^{p,0}(P)(\underbrace{y,\dots,y}_p)=\int_P \inpd{x}{y}^p \dif x, ~~
	M^{0,p}(P)(\underbrace{y,\dots,y}_p)=\int_{\Sn} \ab{\frac{\inpd{x}{y}}{h_P(u)}}^p \dif V_P(u).
\end{align}
Here $\inpd{x}{y}$ denotes the inner product of $x,y$ in the Euclidean space.
These expressions are closely related to the support functions of the $L_p$ moment bodies $M_p P$ and $L_p$ projection bodies $\Pi_p P$:
\begin{align}
h_{M_p P}(y)^p=\int_P \abs{\inpd{x}{y}}^p \dif x, ~~
h_{\Pi_p P}(y)^p=\int_{\Sn} \abs{\frac{\inpd{x}{y}}{h_P(u)}}^p \dif V_P(u).
\end{align}
The $L_p$ moment bodies and $L_p$ projection bodies are fundamental operators in the $L_p$ Brunn-Minkowski Theory (e.g. \cite{HS09a,Lut96,LYZ00a,MR4815903}) and were characterized in the theory of Minkowski valuations; see \cite{Lud02b,Lud05,Hab12b,CLM2017Min,Li2020slnco,LL2016LpMV,Sch2010,SW12,SW2015even,SW2015mink,BL2017Min,HW2025} for related works.
Moreover, as shown by Tang, the first author and Leng \cite{TLL2024} (see also \cite{Li2020slnco}), symmetric-tensor valuations $M^{p,0}, M^{0,p}$ and $L_p$ Minkowski valuations can be uniformly characterized within the framework of function-valued valuations.
That motivates our focus on the asymmetric tensor valuations in this paper.

Let $\MP^n$ be the set of all polytopes in $\R^n$ and $\MPon$ be the subset of polytopes that contain the origin.
The results of Haberl and Parapatits \cite{HP14b,HP16mom,HP17TV} are extended to valuations on $\MPon$ and $\MP^n$ without assuming any regularity by Ludwig and Reitzner \cite{LR2017sl} for $p=0$; by Zeng and the second author \cite{MZ16+Vector}, and the authors and Wang \cite{LMW2020} for $p=1$; and by the second author \cite{Ma21}, and the second author and Wang \cite{MR4230379} for $p=2$.
In this paper, we aim to extend the work of Haberl and Parapatits \cite{HP14b,HP16mom,HP17TV} for $p \ge 2$ to valuations on $\MPon$ and $\MP^n$ removing the regularity assumptions.
More significantly, our classification does not require the symmetry assumptions of tensors and characterizes all non-symmetric tensor valuations of order smaller or equal to $n$.
Note that the results of Haberl and Parapatits \cite{HP14b,HP16mom,HP17TV} cannot follow from our results or from the aforementioned works --- at least we do not know how to do so.  The extension simply means the domain of the valuations is extended from $\MPoon$ to the larger classes $\MPon$ and $\MP^n$ and some additional assumptions (e.g. measurability, symmetry) are removed.

It should be noted that Wang \cite{Wang2022} and Zeng, Zhou \cite{ZZ2024} recently studied the asymmetric case for $p=2$, namely matrix valuations.
They succeeded in giving a representation of matrix valuations determined by decompositions of simplices.
However, they do not have a closed-form expression for matrix valuations for $n=3$ and their representations for $n=2$ are quite complicated.
Consequently, it is hard to extend their results to general tensors; see their representations for $n=2$ in \S \ref{sec:n=2}.

To describe our results, we introduce some definitions and notation.
Every tensor in $\tenset{n}{p}$ can be written as a linear combination of pure tensors of the form $y_1 \otimes \dots \otimes y_p$ with $y_1,\dots,y_p\in \R^n$. For a pure tensor, its action on $(x_1,\dots,x_p) \in (\R^n)^p$ is defined by
\begin{align*}
(y_1 \otimes \dots \otimes y_p) (x_{1}, \dots, x_{p}) = \langle y_1,x_1\rangle \dots  \langle y_p,x_p\rangle.
\end{align*}
Hence we may identify a tensor in $\tenset{n}{p}$ as a $p$-linear functional on $(\R^n)^p$.  
For $\ten{p} \in \tenset{n}{p}$ and $\phi \in \GLn$, we denote by $\phi \cdot \ten{p}$ the natural action of $\phi$ on $\ten{p}$.
That is,
\begin{align*}
	(\phi \cdot \ten{p})(x_{1}, \dots, x_{p})=\ten{p} (\phi^{t}x_{1}, \dots, \phi^{t}x_{p})
\end{align*}
for any $x_1,\dots,x_p \in \R^n$.
A mapping $Z: \mathcal{Q}^n \to \tenset{n}{p}$ is called \emph{$\SLn$ contravariant} if
\begin{align}\label{def:ct1}
Z(\phi P) = \phi^{-t} \cdot ZP
\end{align}
for any $\phi \in \SLn$.

We give a complete classification of $\sln$ contravariant $p$-order-tensor-valued valuations on $\MPon$ for $n \ge p \ge 2$ without assuming that tensors are symmetric.
First, for the case $n-2 \geq p$, only symmetric tensor valuations exist.

\begin{mthm}\label{thm:tensor:p+2}
Let $n-2 \geq p \ge 2$.
A mapping $Z: \MPon \to \tenset{n}{p}$ is an $\SLn$ contravariant valuation if and only if there exists a Cauchy function $\zeta: [0,\infty) \to \R$ such that
\begin{align*}
ZP= \sum_{u \in \mathcal {N}_o (P)} h_P^{-p}(u) \zeta(V(P,u)) \symtp{u}{p}
\end{align*}
for every $P \in \MPon$.
\end{mthm}
Here $\mathcal {N}_o (P)$ is the set of outer unit normals to facets of $P$ which do not contain the origin and $V(P,u):=V_P(\{u\})$. We call $\zeta: [0,\infty) \to \R$ a \emph{Cauchy function} if it satisfies the following Cauchy functional equation:
\begin{align}\label{def:cauchy}
\zeta(r+s)=\zeta(r)+\zeta(s)
\end{align}
for any $r,s \ge 0$. Note that $\zeta(0)=0$. We could write 
$$\sum_{u \in \mathcal {N}_o (P)} h_P^{-p}(u) \zeta(V(P,u)) \symtp{u}{p} = \sum_{u \in S^{n-1} \setminus \{u: h_P(u) = 0\}} h_P^{-p}(u) \zeta(V(P,u)) \symtp{u}{p}.$$

For the case $n-1 =p$, the first non-symmetric tensor valuation appears, which is associated with the cross tensor.
For any $y \in \R^n$, we can define a tensor $y^\times$ of order $n-1$ by
\begin{align}\label{def:cm}
y^\times (x_1,\dots,x_{n-1}) := \det(x_1,\dots,x_{n-1},y).
\end{align}
We call such a tensor $y^\times$ the \emph{cross tensor} of $y$.
The name ``cross" is derived from the case $n=3$, where
\begin{equation}\label{eqn:triple}
	y^\times(x_1,x_2)=\inpd{x_1\times x_2}{y},
\end{equation}
where $x_1 \times x_2$ is the cross product of $x_1$ and $x_2$ in $\R^3$.
We can also use the Hodge star operator on the exterior algebra to define the cross tensor, that is,
$y^\times (x_1,\dots,x_{n-1})=\inpd{y}{\ast(x_1 \wedge \dots \wedge x_{n-1})}$, where $\ast(x_1 \wedge \dots \wedge x_{n-1})$ is the Hodge star operator of $x_1 \wedge \dots \wedge x_{n-1}$ (see \cite[Section 3.3]{Jost2017} for the definition of the Hodge star operator).
For the case $n=2$, $\cross{y}(x)=\inpd{x}{\rho y}$, where $\rho$ is the clockwise rotation in $\R^2$ of the angle $\frac{\pi}{2}$.
Therefore, we can identify $\cross{y} = \rho y$.

Let $m_{n+1}(P)=\int_P x \dif x \in \R^{n}$ be the moment vector of $P$
and $\ct{m}{n+1}{P}$ denote the cross tensor of $m_{n+1}(P)$.

\begin{mthm}\label{thm:tensor:p+1}
Let $n \ge 3$.
A mapping $Z: \MPon \to \tenset{n}{n-1}$ is an $\SLn$ contravariant valuation if and only if there exist a Cauchy function $\zeta: [0,\infty) \to \R$ and a constant $c \in \R$ such that
\begin{align*}
ZP= \sum_{u \in \MN_o(P)}  h_P^{-(n-1)}(u) \zeta(V(P,u)) \symtp{u}{n-1}  + c \ct{m}{n+1}{P}
\end{align*}
for every $P \in \MPon$.
\end{mthm}

For the case $n=p$, non-symmetric tensor valuations also appear. There is an additional class that is linked to the moment vectors of facets instead of the whole body.
The support set of $P$ with outer unit normal vector $u \in S^{n-1}$ is $F(P,u)=\{x\in\R^n: \inpd{x}{u} =h_P(u)\}$ (we allow the dimension of $F(P,u)$ to be less than $n-1$).
Let $m_{n}(F(P,u))=\int_{F(P,u)} x \dif \MH^{n-1}(x)$ be the moment vector of $F(P,u)$, where $\MH^{n-1}$ is the $(n-1)$-dimensional Hausdorff measure.
Denote by $\MN(P)$ the set of all outer unit normals to facets of $P$.
Both the maps
\begin{equation}\label{def:crossfacets}
P \mapsto \sum_{u \in \MN_o(P)} \ct{m}{n}{F(P,u)} \otimes u
\text{~~and~~}
P \mapsto \sum_{u \in \MN(P)} \ct{m}{n}{F(P,u)}  \otimes u
\end{equation}
are $\sln$ contravariant tensor valuations, so are any permutations of them.
Denote by $\symg{p}$ the symmetric group of $\{1,\dots,p\}$.
A permutation $\tau \in \symg{p}$ is a bijection from $\{1,\dots,p\}$ to itself.
We use the two-line notation
\[
\tau = \ab{\begin{smallmatrix}	1  & \cdots & p \\ \tau(1)  & \cdots & \tau(p) \end{smallmatrix}}.
\]
The product of two permutations $\tau_1, \tau_2 \in \symg{p}$ is defined by $\tau_1 \tau_2(i)=\tau_1(\tau_2(i))$ for any $i \in \{1,\dots,p\}$.
For any $\tau \in \symg{p}$, the \emph{permutation} $\tau$ of a $p$-tensor $\ten{p} \in \tenset{n}{p}$ is defined by
\begin{align}\label{def:permtensor}
(\tau \ten{p})(y_1, \dots , y_n):= \ten{p}(y_{\tau^{-1}(1)}, \dots , y_{\tau^{-1}(n)})
\end{align}
for any $y_1,\dots,y_n \in \R^n$.
It is trivial that most permutations of \eqref{def:crossfacets} are linearly dependent.
To classify valuations for the case $p=n$, we have to figure out their nontrivial linear relations, which in turn give some Minkowski-type relations; see Theorem \ref{thm:MR}.

The classification of the case $n=p$ is the following:

\begin{mthm}\label{thm:tensor:p}
Let $n \ge 3$.
A mapping $Z: \MPon \to \tenset{n}{n}$ is an $\SLn$ contravariant valuation if and only if there exist Cauchy functions $\zeta,\eta: [0,\infty) \to \R$ and constants $c_0',c_0, c_1,\dots,c_{n-1}$ such that
\begin{align*}
ZP &= \sum_{u \in \mathcal {N}_o (P)} h_P^{-n}(u) \zeta(V(P,u)) \symtp{u}{n}
+ \sum_{r=0}^{n-2} c_{r+1} \sigma^r \left(\sum_{u \in \MN_o(P)} \ct{m}{n}{F(P,u)} \otimes u \right) \\
&\qquad + \big(c_0'(-1)^{\dim P}V_0(\{o\} \cap \relint P) + c_0V_0(P) + \eta\big(V_n(P)\big) \big) \lct
\end{align*}
for every $P \in \MPon$, where
$\sigma=\ab{\begin{smallmatrix}	1  & \cdots & n-1 & n \\ 2  & \cdots & n & 1 \end{smallmatrix}}$ and $\sigma^r$ is $r$-power of $\sigma$ with respect to the permutation multiplication.

Let $n=2$.
A mapping $Z: \MPo2 \to \tenset{2}{2}$ is an $\SL{2}$ contravariant valuation if and only if there exist Cauchy functions $\zeta,\eta: [0,\infty) \to \R$ and constants $a,b$, $c_0',c_0, c_1$ such that
\begin{align*}
ZP &= \sum_{u \in \mathcal {N}_o (P)} h_P^{-2}(u) \zeta(V(P,u)) \symtp{u}{2}
+ c_{1}  \sum_{u \in \MN_o(P)} \ct{m}{2}{F(P,u) \otimes u }\\
&\qquad + \big(c_0'(-1)^{\dim P}V_0(\{o\} \cap \relint P) + c_0V_0(P) + \eta(V_2(P)) \big) \lct \\
&\qquad + a  M^{2,0}(\rho P)
+ b \dten{2} (\rho P)
\end{align*}
for every $P \in \MPo2$, where $\rho$ is the clockwise rotation in $\R^2$ of the angle $\frac{\pi}{2}$.
\end{mthm}

Here $\lct$ is the Levi-Civita tensor defined in \S \ref{sec:preliminaries}, $V_0$ is the Euler characteristic, and $V_n$ is the $n$-dimensional volume (Lebesgue measure). The valuation $\dten{2}$ only appears in the plane. See definitions of the Euler characteristic and $\dten{2}$ in \S \ref{sec:contra} before Theorem \ref{thm:contrav}.

Remark 1. If we restrict valuations to $\MPoon$, then the valuation
$$P \mapsto \sigma^r \left(\sum_{u \in \MN_o(P)} \ct{m}{n}{F(P,u)} \otimes u \right)$$ would not appear, since the Minkowski-type relations in Theorem \ref{thm:MR} hold. This demonstrates the additional complications that arise when extending the classifications of valuations on $\MPoon$ to $\MPon$.

Remark 2.
The valuations $P \mapsto \ct{m}{n+1}{P}$, $P \mapsto \big(c_0V_0(P) + \eta(V_n(P)) \big) \lct$ and $P \mapsto M^{2,0}(\rho P)$ appeared in Theorems \ref{thm:tensor:p+2}, \ref{thm:tensor:p+1} and \ref{thm:tensor:p} are continuous on $\MPon$ (with respect to the Hausdorff metric on $\MPon$ and the tensor product topology on $\tenset{n}{p}$ induced by the standard topology on $\R^n$). The valuations
$P \mapsto \sum_{u \in \mathcal {N}_o (P)} h_P^{-n}(u) \zeta(V(P,u)) \symtp{u}{n}$, $P \mapsto \left(\sum_{u \in \MN_o(P)} \ct{m}{n}{F(P,u)} \otimes u \right)$ are not continuous on $\MPon$ but continuous on $\MP_{(o)}^n$. For brevity, we omit the details of the proof.

Remark 3.
In light of Henkel and Wannerer \cite{HW2025}, an $\sln$ contravariant valuation $Z: \MPon \to \tenset{n}{p}$ is an invariant valuation under the representation $\Phi: \sln \to \GL(\tenset{n}{p})$, that is, $\Phi(\phi) Z(\phi^{-1} P) = ZP$ for all $P \in \MPon$ and $\phi \in \sln$ with $\Phi(\phi) \ten{p} =\phi^{-t} \cdot \ten{p}$.
Our theorems show that $\sln$ contravariant $\tenset{n}{p}$-valued valuations for $n \ge p \ge 2$ can essentially be decomposed into the corresponding irreducible subrepresentation invariant valuations as $\symtenset{n}{p}$-valued valuations (all $n \ge p$), $\alttenset{n}{p}$-valued valuations ($p=n-1$ or $n$) and $\alttenset{n}{n-1} \otimes \R^n$-valued valuations ($n=p$). Here $\alttenset{n}{p}$ is the space of antisymmetric (alternating) tensors of order $p$ on $\R^n$.
Similar phenomena also appear in Henkel and Wannerer \cite{HW2025}, e.g., \cite[Corollary 1.3]{HW2025} for translation invariant Minkowski valuations.

Finally, let us describe our Minkowski-type relations.

\begin{mthm}\label{thm:MR}
Let $n \ge 2$ and $P \in \MP^n$. We have
\begin{align}\label{eq:MR1}
\sum_{u\in \MN(P)}  \ct{m}{n}{F(P,u)} \otimes u =V_n(P)\lct,
\end{align}
and
\begin{align}\label{eq:MR2}
\sum_{r=1}^n \sgn(\sigma^r)\sigma^r \sum_{u\in \MN(P) \setminus \MN_o(P)}  \ct{m}{n}{F(P,u)} \otimes u=0.
\end{align}
\end{mthm}

By approximation, we can extend \eqref{eq:MR1} to all compact convex sets in $\R^n$ as
\begin{align}
\int_{\partial K} \cross{y} \otimes \nu_K(y) \dif \MH^{n-1}(y) = V_n(K)\lct,
\end{align}
where $\partial K$ is the boundary of $K$, and $\nu_K: \partial K \to \Sn$ is the Gauss map, that is, $\nu_K(y)$ is the outer unit normal to $\partial K$ at $y \in \partial K$ which is well-defined for almost all $y \in \partial K$ with respect to the Hausdorff measure  $\MH^{n-1}$.
Note that $\cross{(y+z)}=\cross{y}+\cross{z}$ for any $y,z \in \R^n$, $V_n(K+z)=V_n(K)$, and $\MH^{n-1}$ is translation invariant.
Hence the previous relation directly implies
$$\int_{\partial K} \cross{z} \otimes \nu_K(y) \dif \MH^{n-1}(y) = \cross{z} \otimes \int_{\Sn}  u \dif S_K(u)=0, ~\forall z \in \R^n,$$
which is equivalent to the classic Minkowski relation $\int_{\Sn} u \dif S_K(u)=0$.
Here $S_K$ is the surface area measure of $K$.

Remark that McMullen \cite{McM97} introduced the following Green-Minkowski connexion:
\begin{equation}
\sum_{u \in \MN(P)} \int_{F(P,u)} x^r \sympro u \dif\MH^{n-1}(x)  = \int_{P} x^{r-1} \sympro \sum_{i=1}^n e_i^2 \dif x
\end{equation}
for $P \in \MP^n$ and $r \ge 1$, where $\sympro$ denotes the symmetric tensor product.
For the case $r=1$ and $n=2$, it is the same as the (tensor) symmetrization of our formula \eqref{eq:MR1}.

The paper is organized as follows. Notation, definitions and results that will be used in this paper are collected in \S \ref{sec:preliminaries}.
In \S \ref{sec:contra}, we prove that the involved operators are $\sln$ contravariant valuations and verify Theorem \ref{thm:MR}.
We prove the planar case $(n=2)$ of Theorem \ref{thm:tensor:p} in \S \ref{sec:n=2} and further prove Theorems \ref{thm:tensor:p+2}, \ref{thm:tensor:p+1} and \ref{thm:tensor:p} for higher dimensional cases $(n>2)$ in \S \ref{sec:n>2}. We extend Theorems \ref{thm:tensor:p+2}, \ref{thm:tensor:p+1} and \ref{thm:tensor:p} to all polytopes not necessarily containing the origin in \S \ref{sec:polyn}.
Finally, in the last section, we briefly list all $\sln$ contravariant tensor valuations with order larger than $n$ that we currently know.

\section{Preliminaries and notation}\label{sec:preliminaries}

Here we collect notation and basics that will be used later; we refer to \cite{Schb2,HP14b,Carr2019,Lud05} for details.

Let $e_1,\dots,e_n$ be the standard basis of $\R^n$. The components of a tensor $\ten{p} \in \tenset{n}{p}$ are denoted by $\ten{p}_{\a_1 \dots \a_p}$ or $\ten{p}(e_{\a_1},\dots,e_{\a_p})$ for $\a_1,\dots,\a_p \in \{1,\dots,n\}$.

We denote by $[A_1,\dots,A_k]$ the convex hull of $A_1,\dots,A_k \subset \R^n$.

Throughout the paper, polytopes and compact convex sets are always assumed to be nonempty.
The support function of a compact convex set $K \subset \R^n$ is
\begin{align}
h_K(u) = \max_{x \in K} \inpd{x}{u}
\end{align}
for any $u \in \R^n$.
Usually, the cone-volume measure is defined for a compact convex set $K$ containing the origin.
Here we need to extend its definition to all compact convex sets so that we can extend our results to $\MP^n$.
First, the surface area measure of a full-dimensional compact convex set $K$ is defined by $S_K(\omega):=\MH^{n-1} (\nu_K^{*}(\omega))$ for any Borel set $\omega \subset \Sn$ where $\nu_K^{*}(\omega)=\{x \in \partial K: \inpd{x}{u} = h_K(u) \text{~for~some~} u \in \omega\}$ is the reverse Gauss image of $\omega$.
For a lower-dimensional compact convex set $K$, its surface area measure is defined by $S_K=\MH^{n-1}(K) (\delta_{-u_0} + \delta_{u_0})$ where $u_0 \in S^{n-1}$ is orthogonal to the affine hull of $K$ and $\delta_{u_0}$ is the Dirac measure at $u_0$.
Then the cone-volume measure of a compact convex set $K$ is a signed measure defined by $V_K(\omega):=\frac{1}{n}\int_{\omega} h_K(u) \dif S_K(u)$.

We also call $\zeta: \R \to \R$ a \emph{Cauchy function} if it satisfies the following Cauchy functional equation:
\begin{align}
\zeta(r+s)=\zeta(r)+\zeta(s)
\end{align}
for any $r,s \in \R$.
Clearly, if $\zeta$ is a Cauchy function on $[0,\infty)$, then its unique extension to a Cauchy function on $\R$ is defined by $\zeta(-r):=-\zeta(r)$ for $r>0$.
Note that a Cauchy function on $\R$ is not necessarily linear, but it is linear if we assume some regularity conditions, e.g. measurability, continuity, monotonicity, boundedness, etc.

\subsection{The decomposition of $\sln$ contravariant valuations}\label{sec:DV}

Let $\slpm{n}=\{\phi \in \gln : \det \phi = 1 \text{~or~} -1\}$.
For $\delta \in \{0,1\}$, we say $Z:\mathcal{Q}^n \to \tenset{n}{p}$ is $\slpm{n}$-$\delta$-contravariant if
\begin{align*}
Z(\phi P) =(\det\phi)^\delta (\phi^{-t} \cdot ZP)
\end{align*}
for every $P\in \mathcal Q^n$ and $\phi \in \slpm{n}$.

Denote by $\tenval(\mathcal{Q}^n;\tenset{n}{p})$ the set of $\sln$ contravariant valuations $Z:\mathcal{Q}^n \to \tenset{n}{p}$ and by $\tenval_\delta(\mathcal{Q}^n;\tenset{n}{p})$ the set of $\slpm{n}$-$\delta$-contravariant valuations $Z:\mathcal{Q}^n \to \tenset{n}{p}$.
Choose a $\theta \in \slpm{n}$ with $\det \theta=-1$.
We can decompose a valuation $Z \in \tenval(\mathcal{Q}^n;\tenset{n}{p})$ as
$Z=Z^+ + Z^-,$
where
\begin{align*}
&Z^+P:=\frac{1}{2}\ab{ZP+\theta^t \cdot Z(\theta P)},\\
&Z^-P:=\frac{1}{2}\ab{ZP-\theta^t \cdot Z(\theta P)}.
\end{align*}
Since $Z$ is an $\sln$ contravariant valuation, 
we easily find that $Z^+$ and $Z^-$ do not depend on the choice of $\theta$, and $Z^+ \in \tenval_0(\mathcal{Q}^n;\tenset{n}{p})$, $Z^- \in \tenval_1(\mathcal{Q}^n;\tenset{n}{p})$, i.e.,
\begin{align*}
	\tenval(\mathcal{Q}^n;\tenset{n}{p}) = \tenval_0(\mathcal{Q}^n;\tenset{n}{p}) \oplus \tenval_1(\mathcal{Q}^n;\tenset{n}{p}).
\end{align*}
The idea of this decomposition is due to Haberl and Parapatits \cite{HP17TV}.

\subsection{The Levi-Civita tensor}
Let $\sgn$ denote the sign of a permutation.
The Levi-Civita symbols $\lct_{\a_1 \dots \a_{n}}$ on $\R^n$
for $\a_1, \dots ,\a_{n} \in \{1,\dots,n\}$ are defined as follows:
\begin{align}\label{def:LV}
\lct_{\a_1 \dots \a_{n}}
:=\begin{cases}
\sgn \ab{\begin{smallmatrix}
1  & \cdots & n \\
\a_1  & \cdots & \a_n
\end{smallmatrix}}, & \text{if~} \a_1,\dots,\a_{n}~\text{are~pairwise~distinct},\\
0,  & \text{otherwise}.
\end{cases}
\end{align}
The Levi-Civita tensor $\lct$ on $\R^n$ is the $n$-order tensor whose components are Levi-Civita symbols.
Note that
$\lct_{\a_1 \dots \a_{n}}= \det (e_{\a_1},\dots,e_{\a_n})$,
we have
\begin{align}
\lct (x_1,\dots,x_n)= \det(x_1,\dots,x_n)
\end{align}
for any $x_1,\dots,x_n \in \R^n$.
Thus
\begin{align}
\phi^{-t} \cdot \lct (x_1,\dots,x_n)
&= \lct (\phi^{-1} x_1,\dots,\phi^{-1} x_n)
= \det(\phi^{-1} x_1,\dots,\phi^{-1} x_n) \\
&= (\det\phi)^{-1} \det(x_1,\dots,x_n) = (\det\phi)^{-1} \lct (x_1,\dots,x_n)
\end{align}
for any $\phi \in \slpm{n}$.
This shows that the mapping $P \mapsto \lct$ for all $P \in \MPon$ is not only $\slpm{n}$ invariant, but also $\slpm{n}$-$1$-contravariant.

In this paper, for example, in \eqref{eq615-1} and \eqref{eq615-2}, we need to frequently involve the Levi-Civita symbol $\lct_{\a_1,\dots,\a_n}$, where $\a_1,\dots,\a_{n-1}$ are specified in advance, and $\a_n \in \{1,\dots,n\} \setminus \{\a_1,\dots,\a_{n-1}\}$. This forces us to frequently emphasize what $\a_n$ is, which is rather verbose. Therefore, we extend the Levi-Civita symbol and define the extended Levi-Civita symbol as follows:
\[
\lct_{\a_1 \dots \a_{n-1}}:=\lct_{\a_1 \dots \a_{n-1} \a_{n}},
\]
where $\a_n \in \{1,\dots,n\} \setminus \{\a_1,\dots,\a_{n-1}\}$. Since the ambient space is $\mathbb{R}^n$ by default, this notation should not cause confusion.

For a given $i \in \seq{1}{n}$, we will also need to consider the Levi-Civita symbols on the space $e_i^\bot:= \{x \in \R^n : \inpd{x}{e_i} = 0\} \cong \mathbb{R}^{n-1}$, which are defined by
\begin{align*}
\lct_{\a_1 \dots \a_{n-1}}^{e_i^\bot}:=\det\nolimits_{n-1} (e_{\a_1},\dots,e_{\a_{n-1}})
\end{align*}
for any $\a_1, \dots, \a_{n-1} \in \seq{1}{n} \setminus \set{i}$, and similarly we can define the extended Levi-Civita symbols on $e_i^\bot$ by
\begin{align*}
\lct_{\a_1 \dots \a_{n-2}}^{e_i^\bot}:=\lct_{\a_1 \dots \a_{n-1}}^{e_i^\bot}
\end{align*}
for any $\a_1, \dots, \a_{n-2} \in \seq{1}{n} \setminus \set{i}$ and $\a_{n-1} \in \seq{1}{n} \setminus \{\seqalpha{\a}{1}{n-2},i\}$.

\subsection{The dissections of the standard simplex}

We denote the standard simplex by $T^{d}=[o,e_1,\dots,e_{d}]$ for any $1 \leq d \leq n$.
Let $\hat{e}_{i}$ denote that $e_{i}$ is omitted and let 
$\hat T^{d-1}_i=[o,e_1,\dots,\hat{e}_{i},\dots,e_d]$ for $1 \leq i \le d$.

The following triangulations are used several times in the proofs.
For $0 < \lambda <1$, the hyperplane $H_\lambda$ and half spaces $H_\lambda^-,H_\lambda^+$ are defined by
\begin{align*}
H_\lambda := \{x \in\R^n : \inpd{x}{((1-\lambda) e_1- \lambda e_2)} = 0\}, \\
H_\lambda^- := \{ x \in\R^n : \inpd{x}{((1-\lambda) e_1- \lambda e_2)} \leq 0 \}, \\
H_\lambda^+ := \{ x \in\R^n : \inpd{x}{((1-\lambda) e_1- \lambda e_2)} \geq 0 \}.
\end{align*}
Let $Z\in \tenval(\MPon;\tenset{n}{p})$. We have
\begin{align}\label{val}
Z (sT^{d}) + Z (sT^{d} \cap H_\lambda) = Z (sT^{d} \cap H_\lambda ^-) + Z (sT^{d} \cap H_\lambda ^+)
\end{align}
for $s > 0$ and $2 \leq d \leq n$.
If $Z$ is simple, then $Z (sT^{d} \cap H_\lambda) =0$.

For the case $d <n$, we define $\phi_1,\psi_1 \in \sln$ by
\begin{align*}
\phi_1:=\ab{\begin{array}{cccccc}
\lambda & & & & & \\
1-\lambda &1 & & & &\\
&& 1 && &\\
&&&\ddots & &\\
&&& & 1 &\\
&&&&& \frac{1}{\lambda}
\end{array}}, \qquad
\psi_1:=\ab{\begin{array}{cccccc}
1 &\lambda\ & & & & \\
 &1-\lambda & & & &\\
&& 1 && &\\
&&&\ddots &&\\
&&& & 1 &\\
&&&&& \frac{1}{1-\lambda}
\end{array}}.
\end{align*}
Notice that $T^d\cap H_\lambda ^- = \phi_1 T^d$, $T^d\cap H_\lambda ^+ = \psi_1 T^d$ and $T^d \cap H_\lambda =\phi_1 \hat T^{d-1}_2$.
Applying the $\SLn$ contravariance of $Z$ in \eqref{val}, we have
\begin{equation}\label{30a}
\begin{aligned}
&Z (sT^d)(e_{\a_1}, \dots, e_{\a_p}) + Z (s\hat T^{d-1}_2)(\phi_1^{-1}e_{\a_1}, \dots, \phi_1^{-1}e_{\a_p}) \\
&=Z (sT^d)(\phi_1^{-1} e_{\a_1}, \dots, \phi_1^{-1}e_{\a_p})+ Z(sT^d)(\psi_1^{-1} e_{\a_1}, \dots, \psi_1^{-1}e_{\a_p})
\end{aligned}
\end{equation}
for $2 \le d \le n-1$, where
\begin{align*}
\phi_1^{-1}=\ab{\begin{array}{cccccc}
\frac{1}{\lambda} & & & & &\\
-\frac{1-\lambda}{\lambda} &1 & & &&\\
&& 1 &&& \\
&&&\ddots &&\\
&&&&1&\\
&&&&&\lambda
\end{array}}, \qquad
\psi_1^{-1}=\ab{\begin{array}{cccccc}
1 &-\frac{\lambda}{1-\lambda} & & && \\
 &\frac{1}{1-\lambda} & & &&\\
&& 1 &&& \\
&&&\ddots &&\\
&&&&1&\\
&&&&&1-\lambda
\end{array}}.
\end{align*}

For the full-dimensional case, set $\phi_2,\psi_2 \in \gln$ as
\begin{align*}
\phi_2:=\ab{\begin{array}{ccccc}
\lambda & & & & \\
1-\lambda &1 & & &\\
&& 1 && \\
&&&\ddots &\\
&&&&1
\end{array}}, \qquad
\psi_2:=\ab{\begin{array}{ccccc}
1 &\lambda\ & & & \\
 &1-\lambda & & &\\
&& 1 && \\
&&&\ddots &\\
&&&&1
\end{array}}.
\end{align*}
Assume that $Z$ is an $\SLn$ contravariant valuation. Then
\begin{equation}\label{30}
\begin{aligned}
&Z (sT^{n})_{\a_1 \dots \a_p} + \lambda^{p/n} Z (\lambda^{1/n}s\hat T^{d-1}_2)(\phi_1^{-1}e_{\a_1}, \dots, \phi_1^{-1}e_{\a_p})\\
&= Z (\phi_2 sT^{n})_{\a_1 \dots \a_p}+ Z(\psi_2 sT^n)_{\a_1 \dots \a_p} \\
&=\lambda^{p/n} Z (\lambda^{1/n} sT^{n})(\phi_2^{-1} e_{\a_1}, \dots, \phi_2^{-1}e_{\a_p})\\
&\qquad + (1-\lambda)^{p/n}Z((1-\lambda)^{1/n} sT^n)(\psi_2^{-1} e_{\a_1}, \dots, \psi_2^{-1}e_{\a_p}),	
\end{aligned}
\end{equation}
for every $s>0$, where
\begin{align*}
\phi_2^{-1}=\ab{\begin{array}{ccccc}
\frac{1}{\lambda} & & & & \\
-\frac{1-\lambda}{\lambda} &1 & & &\\
&& 1 && \\
&&&\ddots &\\
&&&&1
\end{array}}, \qquad
\psi_2^{-1}=\ab{\begin{array}{ccccc}
1 &-\frac{\lambda}{1-\lambda} & & & \\
 &\frac{1}{1-\lambda} & & &\\
&& 1 && \\
&&&\ddots &\\
&&&&1
\end{array}}.
\end{align*}

The following lemmas tell us that $\SLn$ contravariant tensor valuations are uniquely determined by their values on simplices. The proofs are similar to the brief argument of \cite[Lemma 4.5 and Lemma 4.6]{Li2018AFV}, or see \cite[\S 2 and \S 4]{LR2017sl} for a similar detailed argument.
\begin{lem}\label{lemuq}
Let $Z$ and $Z'$ be $\SLn$ contravariant tensor valuations on $\MPon$. If $Z (sT^d) = Z' (sT^d)$ for every $s>0$ and $0 \leq d \leq n$, then $Z P = Z' P$ for every $P \in \MPon$.
\end{lem}

\begin{lem}\label{lemuq2}
Let $Z$ and $Z'$ be $\SLn$ contravariant tensor valuations on $\MP^n$. If $Z (sT^d) = Z' (sT^d)$ and $Z(s[e_1,\dots,e_d]) = Z'(s[e_1,\dots,e_d])$ for every $s>0$ and $0 \leq d \leq n$, then $Z P = Z' P$ for every $P \in \MP^n$.
\end{lem}

\subsection{Classifications of $\sln$ contravariant tensor valuations for orders $0$ and $1$.}
In the proof, we require the following classifications of $\sln$ contravariant tensor valuations of orders $0$ (real-valued valuations) and $1$ (vector-valued valuations) established by Ludwig and Reitzner \cite{LR2017sl}, and the authors and Wang \cite{LMW2020}, respectively.
\begin{thm}[\cite{LR2017sl}]\label{thm:p=0}
Let $n \ge 2$. A mapping $Z: \MPon \to \R$ is an $\SLn$ invariant valuation if and only if
there exist constants $c_0',c_0\in\R$ and a Cauchy function $\eta: [0,\infty) \to \R$ such that
\begin{align*}
ZP=c_0'(-1)^{\dim P}V_0(\{o\} \cap \relint P) + c_0V_0(P) + \eta\big(V_n(P)\big),
\end{align*}
for every $P \in \MPon$.
\end{thm}

\begin{thm}[\cite{LMW2020}]\label{thm:p=1}
Let $n \ge 3$. A mapping $Z: \MPon \to \R^n$ is an $\SLn$ contravariant valuation if and only if there exists a Cauchy function $\zeta: [0,\infty) \to \R$ such that
\begin{align*}
ZP= \sum_{u \in \mathcal {N}_o (P)} h_P^{-1}(u)\zeta(V(P,u)) u
\end{align*}
for every $P \in \MPon$.

A mapping $ Z :\polyeo\to\R^2$ is an $\SL{2}$ contravariant valuation
if and only if there exist constants $b,c\in\R$
and a Cauchy function $\zeta: [0,\infty) \to \R$ such that
\[ ZP=\sum_{u \in \mathcal {N}_o (P)} h_P^{-1}(u)\zeta(V(P,u)) u
+ c \rho m_3(P) + b\dten{1}(\rho P)\]
for every $P\in\polyeo$.
\end{thm}

Here $\dten{1}(\rho P)$ is a vector valuation defined similarly to the matrix valuation $\dten{2}(\rho P)$; see the definitions of such tensor (not merely vector and matrix) valuations in \S \ref{sec:contra} before Theorem \ref{thm:contrav}.

Remark:
Let $P \in \MP^n$. As mentioned in the first section, we have $\ct{m}{3}{P}=\rho m_3(P)$ in the plane.
Although it is an easy result, it is an important observation that allows us to proceed with induction to classify tensor valuations.
Another observation is $\zeta(V_n(P))=\sum_{u \in \mathcal {N}_o (P)} \zeta(V(P,u))$. Hence, for $0 \le p < \infty$, the sum $\sum_{u \in \mathcal {N}_o (P)} h_P^{-p}(u)\zeta(V(P,u)) u^p$ is the tensor generalization of the scalar $\zeta(V_n(P))$ and the vector $\sum_{u \in \mathcal {N}_o (P)} h_P^{-1}(u)\zeta(V(P,u)) u$, reducing to them when $p=0$ and $p=1$ respectively.

\section{Valuations and contravariances}\label{sec:contra}
First, note that the moment vector is a valuation which is \emph{$\slpm{n}$ covariant}, that is
\begin{align}
m_{n+1}(\phi P) = \phi m_{n+1}(P)
\end{align}
for any $\phi \in \slpm{n}$ and $P \in \MP^n$.
Hence the mapping $P \mapsto \ct{m}{n+1}{P}$ is in $\tenval_1(\MP^n;\tenset{n}{p})$. In fact,
\begin{align}
\big(\ct{m}{n+1}{\phi P} \big) (x_1,\dots,x_{n-1}) = \det(x_1,\dots,x_{n-1},m_{n+1}(\phi P))= (\det \phi) \ab{\phi^{-t} \cdot \ct{m}{n+1}{P}}.
\end{align}

For $P \in \mpn$, let $\MN(P)$ be the set of outer unit normals of $P$ and $\MN_o(P)$ be the set of outer unit normals of $P$ such that the affine hulls of the corresponding facets do not contain the origin. We have the following.

\begin{lem}\label{lem:mo1}
For any $P \in \MP^n$, $u \in \MN(P)$ and $\phi \in \slpm{n}$, we have
\begin{align*}
m_n\left(F\left(\phi P,\overline{\phi^{-t} u}\right)\right)=| \phi^{-t} u|  \phi m_n\left(F\left(P,u\right)\right),
\end{align*}
where $\overline{\phi^{-t} u}= \frac{\phi^{-t} u}{|\phi^{-t} u|}$.
\end{lem}
\begin{proof}
It is clear that, for any $P \in \MP^n$, $u \in \MN_o(P)$ and $\phi \in \slpm{n}$,
\begin{equation}\label{eqn:hf}
	h_{\phi P}(\overline{\phi^{-t}u})=|\phi^{-t}u|^{-1}h_P(u)\quad\text{and}\quad F\left(\phi P,\overline{\phi^{-t}u}\right)=\phi F(P,u).
\end{equation}
Since $m_n(\lambda K)= \lambda^{n} m_n(K)$ for any $\lambda >0$ and convex body $K \subset \R^n$ with $\dim K =n-1$, we use Fubini's theorem to obtain
\begin{equation}\label{eqn:mvof}
	m_{n+1}([o,F(P,u)])=\frac1{n+1}h_P(u)m_n(F(P,u)),
\end{equation}
which implies
\begin{equation}
\begin{split}
	m_{n+1}(\phi[o,F\left(P,u\right)])
	&=m_{n+1}([o,F\left(\phi P,\overline{\phi^{-t}u}\right)]) \\
	&=	\frac1{n+1}h_{\phi P}(\overline{\phi^{-t}u})m_n\left(F\left(\phi P,\overline{\phi^{-t}u}\right)\right)\\
	&= \frac1{n+1}|\phi^{-t}u|^{-1}h_P(u)m_n\left(F\left(\phi P,\overline{\phi^{-t}u}\right)\right),
\end{split}
\end{equation}
due to \eqref{eqn:hf}.
On the other hand,
\begin{equation}
	m_{n+1}\left(\phi[o,F(P,u)]\right) = \phi m_{n+1}([o,F(P,u)]).
\end{equation}
Therefore, the desired result for $u \in \MN_o(P)$ follows.

If $u \in \MN(P) \setminus \MN_o(P)$, we can translate $P$ by some suitable sequence $(x_k)_{k=1}^\infty \in \R^n$ such that $u \in \MN_o(P+x_k)$ for all $k$ and $x_k \to o$ as $k \to \infty$.
Then by the result for $u \in \MN_o(P+x_k)$, we have
\begin{align*}
m_n\left(F\left(\phi (P+x_k),\overline{\phi^{-t} u}\right)\right)=| \phi^{-t} u|  \phi m_n\left(F\left(P+x_k,u\right)\right),
\end{align*}
Since
\begin{align}
\lim_{k \to \infty} m_n\left(F\left(\phi (P+x_k),\overline{\phi^{-t} u}\right)\right)=m_n\left(F\left(\phi P,\overline{\phi^{-t} u}\right)\right)
\end{align}
and
\begin{align}
\lim_{k \to \infty} m_n\left(F\left(P+x_k,u\right)\right) =m_n\left(F\left(P,u\right)\right),
\end{align}
the desired result hold for $u \in \MN(P) \setminus \MN_o(P)$.
\end{proof}

Lemma \ref{lem:mo1} implies contravariance of the following valuations.
\begin{lem}\label{lem:mo2}
The mappings
\begin{equation}
\MP^n \ni P \mapsto \sum_{u \in \MN_o(P)} \ct{m}{n}{F(P,u)} \otimes u,
\text{~~and~~}
\MP^n \ni P \mapsto \sum_{u \in \MN(P)} \ct{m}{n}{F(P,u)}  \otimes u
\end{equation}
are in $\tenval_1(\MP^n;\tenset{n}{n})$.
\end{lem}
\begin{proof}
Let $P, Q \in \MP^n$ such that $P \cup Q \in \MP^n$.
We can divide $S^{n-1}$ into three parts: $\omega_1=\{u \in S^{n-1}: h_P(u) < h_Q(u)\}$, $\omega_2=\{u \in S^{n-1}: h_P(u) > h_Q(u)\}$ and $\omega_3=\{u \in S^{n-1}: h_P(u) = h_Q(u)\}$.
For any $u \in S^{n-1}$, observe that
\begin{align}
h_{P \cap Q}(u) = \min\{h_P(u),h_Q(u)\} \text{~and~} h_{P \cup Q}(u) = \max\{h_P(u),h_Q(u)\},
\end{align}
\begin{align}
F(P \cap Q,u) =
\begin{cases}
F(P,u), &\text{~if~} u \in \omega_1, \\
F(Q,u), &\text{~if~} u \in \omega_2, \\
F(P,u) \cap F(Q,u), &\text{~if~} u \in \omega_3,
\end{cases}
\end{align}
and
\begin{align}
F(P \cup Q,u) =
\begin{cases}
F(Q,u), &\text{~if~} u \in \omega_1, \\
F(P,u), &\text{~if~} u \in \omega_2, \\
F(P,u) \cup F(Q,u), &\text{~if~} u \in \omega_3.
\end{cases}
\end{align}
Hence
\begin{align}
\MN_o(P \cap Q) \cap \omega_1=\MN_o(P) \cap \omega_1, ~\MN_o(P \cap Q) \cap \omega_2=\MN_o(Q) \cap \omega_2, \\
\MN_o(P \cup Q) \cap \omega_1=\MN_o(Q) \cap \omega_1, ~\MN_o(P \cup Q) \cap \omega_2=\MN_o(P) \cap \omega_2.
\end{align}
and therefore
\begin{align}
&\sum_{u \in \MN_o(P\cap Q) \cap (\omega_1 \cup \omega_2)} \ct{m}{n}{F(P\cap Q,u)} \otimes u + \sum_{u \in \MN_o(P\cup Q) \cap (\omega_1 \cup \omega_2)} \ct{m}{n}{F(P\cup Q,u)} \otimes u\\
&=\sum_{u \in \MN_o(P) \cap (\omega_1 \cup \omega_2)} \ct{m}{n}{F(P,u)} \otimes u + \sum_{u \in \MN_o(Q) \cap (\omega_1 \cup \omega_2)} \ct{m}{n}{F(Q,u)} \otimes u.
\end{align}
For $u \in \omega_3$, some $F(P,u)$ or $F(Q,u)$ might be less than $(n-1)$-dimensional. Note that $m_n(K\cap L)+m_n(K\cup L)= m_n(K)+m_n(L)$ for any $K,L \in \MP^n$ with $\dim K = \dim L \le n-1$ and $m_n$ vanishes for any convex set with dimension less than $n-1$, so is the corresponding cross tensor.
Therefore, we still have
\begin{align}
&\sum_{u \in \MN_o(P\cap Q) \cap \omega_3} \ct{m}{n}{F(P\cap Q,u)} \otimes u + \sum_{u \in \MN_o(P\cup Q) \cap \omega_3} \ct{m}{n}{F(P\cup Q,u)} \otimes u\\
&=\sum_{u \in \MN_o(P) \cap \omega_3} \ct{m}{n}{F(P,u)} \otimes u + \sum_{u \in \MN_o(Q) \cap \omega_3} \ct{m}{n}{F(Q,u)} \otimes u.
\end{align}
Thus
$P \mapsto \sum_{u \in \MN_o(P)} \ct{m}{n}{F(P,u)} \otimes u$ is a tensor valuation.

By Lemma \ref{lem:mo1} and the fact that $u \in \MN_o(P)$ if and only if $v=\overline{\phi^{-t} u}\in \MN_o(\phi P)$ for any $\phi \in \slpm{n}$, we have
\begin{align*}
&\ab{\sum_{v \in \MN_o( \phi P)} \ct{m}{n}{F(\phi P,v)} \otimes v}(x_1,\dots,x_{n-1},x_n) \\
&=\ab{\sum_{u \in \MN_o(P)} \ct{m}{n}{F(\phi P,\overline{\phi^{-t} u})} \otimes \overline{\phi^{-t} u} }(x_1,\dots,x_{n-1},x_n)\\
&=\sum_{u \in \MN_o(P)} \det \ab{x_1,\dots,x_{n-1},m_n(F(\phi P,\overline{\phi^{-t} u}))} \inpd{x_n}{\overline{\phi^{-t} u}}\\
&= \sum_{u \in \MN_o(P)} \det(\phi)  \det \ab{\phi^{-1}x_1,\dots,\phi^{-1}x_{n-1},m_n(F(P, u))} \inpd{\phi^{-1}x_n}{u}\\
&= (\det \phi) \ab{\phi^{-t} \cdot \sum_{u \in \MN_o(P)} \ab{\ct{m}{n}{F(P,u)} \otimes u}}.
\end{align*}
Thus the mapping $P \mapsto \sum_{u \in \MN_o(P)} \ct{m}{n}{F(P,u)} \otimes u$ is in $\tenval_1(\MP^n;\tenset{n}{p})$.

Similarly, the mapping $P \mapsto \sum_{u \in \MN(P)} \ct{m}{n}{F(P,u)} \otimes u$ is in $\tenval_1(\MP^n;\tenset{n}{p})$.
\end{proof}

The Euler characteristic $V_n(P)=1$ for nonempty $P \in \MP^n$ and $V_n(\varnothing)=0$.
Hence it is clear an $\SLn$ invariant valuation.
The mapping $\dten{p}$ for $p \ge 1$ in Theorem \ref{thm:tensor:p} is defined for $P \in \MP^2$ by
$$\dten{p}(P):=v^p+w^p,$$
if $\dim P=2$ and $P$ has two edges $[o,v]$ and $[o,w]$,
or $\dim P=2$ and $P$ has an edge $[v,w]$ such that $o \in (v,w)$;
$$\dten{p}(P):=2(v^p+w^p),$$
if $\dim P=1$ and $P=[v,w]$ contains the origin; and
$$\dten{p}(P):=0,$$ otherwise.
By an argument similar to that in the proof of \cite[Lemma 2.4]{MZ16+Vector}, we can verify that it is a valuation. 

Recall that $\rho$ is the clockwise rotation in $\R^2$ of the angle $\frac{\pi}{2}$. We have $\rho \phi \rho^{-1} = \det(\phi) \phi^{-t}$ for every $\phi \in \slpm{2}$.
Now we summarize the contravariances of the valuations to be classified. The proof of the valuation property of the mapping $\MP^n \ni P \mapsto \sum_{u \in \mathcal {N}_o (P)}  \zeta(V(P,u)) \symtp{\ab{\frac{u}{h_P(u)}}}{p}$ can be found in \cite[Corollary 3 and Lemma 3]{TLL2024}. 
The further proofs are similar to the proof of Lemma \ref{lem:mo2}.
\begin{thm}\label{thm:contrav}
Let $n \ge 2$, $\zeta,\eta: [0,\infty) \to \R$ be Cauchy functions, and $c_0',c_0 \in \R$ be constants.
We have the following.
\begin{align}
& \abb{\MP^n \ni P \mapsto \sum_{u \in \mathcal {N}_o (P)}  \zeta(V(P,u)) \symtp{\ab{\frac{u}{h_P(u)}}}{p} } \in  \tenval_0(\MP^n;\tenset{n}{p}), \\
& \abb{\MP^2 \ni P \mapsto  M^{2,0}(\rho P)} \in \tenval_0(\MP^2;\tenset{2}{2}) ,\\
& \abb{\MP^2 \ni P \mapsto \dten{2} (\rho P)} \in \tenval_0(\MP^2;\tenset{2}{2}), \\
& \abb{\MP^n \ni P \mapsto c_0'(-1)^{\dim P}V_0(\{o\} \cap \relint P) + c_0V_0(P) + \eta(V_n(P))} \in \tenval_0(\MP^n;\R),
\end{align}
and
\begin{align}
&\abb{\MP^n \ni P \mapsto\ct{m}{n+1}{P}} \in \tenval_1(\MP^n;\tenset{n}{n-1})\\
&\abb{\MP^n \ni P \mapsto \sum_{u \in \MN_o(P)} \ct{m}{n}{F(P,u)} \otimes u} \in \tenval_1(\MP^n;\tenset{n}{n})\\
& \abb{\MP^n \ni P \mapsto \sum_{u \in \MN(P)} \ct{m}{n}{F(P,u)}  \otimes u}  \in \tenval_1(\MP^n;\tenset{n}{n})\\
& \abb{\MP^n \ni P \mapsto \big(c_0'(-1)^{\dim P}V_0(\{o\} \cap \relint P) + c_0V_0(P) + \eta(V_n(P))\big) \lct} \in \tenval_1(\MP^n;\tenset{n}{n})\\
& \abb{\MP^n \ni P \mapsto \dten{1}(\rho P)} \in \tenval_1(\MP^n;\R^n).
\end{align}
\end{thm}

Next, we calculate the moment vectors of the standard simplex and its faces. Denote $\bar e := \frac{1}{\sqrt n} \ab{e_1+\dots+e_n}$.

\begin{lem}\label{lem:mo3}
Let $d \le n$.
Define $m_{d+1}(T^{d}) := \int_{T^d} x \dif \MH^{d}(x)$. Then
\begin{align*}
m_{d+1}(T^{d}) =\frac{1}{(d+1)!}(e_1+\dots+e_d).
\end{align*}
Moreover, we have
\begin{align*}
m_n(F(T^n,-e_{i}))=\frac{1}{n!} (e_1+\dots+\hat{e}_{i}+\dots+e_n),
\end{align*}
for any $i \in \{1,\dots,n\}$ and
\begin{align*}
m_n(F(T^n,\bar e))=\frac{\sqrt{n}}{n!} (e_1+\dots+e_n).
\end{align*}
\end{lem}
\begin{proof}
For $x \in \R^n$, we write $x=\sum_{i=1}^n x^i e_i$ for $x^i \in \R$ and $1 \le i \le n$.
By Fubini's theorem, for any $1 \leq d \leq n$,
\begin{align*}
m_{d+1}(T^{d})= \int_{T^d} x \dif x &=\sum_{i=1}^d \left(\int_{T^d} x^i \dif x \right) e_i\\
&=\sum_{i=1}^d \left(\int_0^1 \frac{1}{(d-1)!}s(1-s)^{d-1} \dif s \right)e_i \\
&=\frac{1}{(d+1)!}(e_1+\dots+e_d).
\end{align*}
Hence we obtain the first formula and then the second formula follows from $F(T^n,-e_{i})=[o,e_1,\dots,\hat{e}_{i},\dots,e_n]$.
By \eqref{eqn:mvof}, we get
\begin{align*}
m_n(F(T^n,\bar e))&=\frac{(n+1)m_n(T^n)}{h_{T^n}(\bar e)}=\frac{\sqrt{n}}{n!} (e_1+\dots+e_n).
\end{align*}
\end{proof}

We can now prove Theorem \ref{thm:MR}.
\begin{proof}[Proof of Theorem \ref{thm:MR}]
The formula \eqref{eq:MR1} can be proved directly with the divergence theorem.
Here we use the $\sln$ contravariance and valuation property to prove both \eqref{eq:MR1} and \eqref{eq:MR2}.
Since both sides of \eqref{eq:MR1} and \eqref{eq:MR2} are simple, $n$-homogeneous ($Z(\lambda P) = \lambda^n ZP$ for any $\lambda >0$) and $\sln$ contravariant valuations, by Lemma \ref{lemuq2}, we only need to show that \eqref{eq:MR1} and \eqref{eq:MR2} hold for the standard simplex $T^n$.
We first prove \eqref{eq:MR2} and then \eqref{eq:MR1}.

1. We show that \eqref{eq:MR2} holds for the standard simplex $T^n$.
By Lemma \ref{lem:mo3}, we have
\begin{align}\label{eq:1}
&\left(\sum_{u\in \MN(T^n) \setminus \MN_o(T^n)}  \ct{m}{n}{F(T^n,u)} \otimes u \right)  (e_{\a_1}, \dots, e_{\a_n}) \notag\\
&=\sum_{i=1}^n  (-e_i \cdot e_{\a_n}) \det (e_{\a_1},\dots,e_{\a_{n-1}}, m_n(F(T^n,-e_i))) \notag\\
&=- \det (e_{\a_1},\dots,e_{\a_{n-1}}, m_n(F(T^n,-e_{\a_n}))) \notag\\
&=- \frac{1}{n!} \det (e_{\a_1},\dots,e_{\a_{n-1}},(e_1+\dots+\hat{e}_{\a_n}+\dots+e_n))
\end{align}
for any $\a_1,\dots,\a_n \in \{1,\dots,n\}$.
Hence
$$\left(\sum_{u\in \MN(T^n) \setminus \MN_o(T^n)}  \ct{m}{n}{F(T^n,u)} \otimes u \right)   (e_{\a_1}, \dots, e_{\a_n}) \neq 0$$
if only if there exists only one $i \in \{1,\dots,n-1\}$ such that $\a_i=\a_n$.
Note that for any $\ten{p} \in \tenset{n}{p}$, $\sigma^r \ten{p}=0$ if and only if $\ten{p}=0$.
To obtain \eqref{eq:MR2}, it is now sufficient to show that
\begin{align*}
0&=\left(\sum_{u\in \MN(T^n) \setminus \MN_o(T^n)}  \ct{m}{n}{F(T^n,u)} \otimes u \right)  (e_{\a_1}, \dots, e_{\a_n}) \\
&\qquad +\sgn(\sigma^{n-i})\left(\sum_{u\in \MN(T^n) \setminus \MN_o(T^n)}  \ct{m}{n}{F(T^n,u)} \otimes u \right) \ab{e_{\a_{\sigma^{-(n-i)}(1)}}, \dots, e_{\a_{\sigma^{-(n-i)}(n)}}},
\end{align*}
if $\a_i=\a_n$ for some $i \in \{1,\dots,n-1\}$ since \eqref{def:permtensor} and $\sigma^0=\sigma^n$.
Using \eqref{eq:1} again, we only need to show that
\begin{align*}
&\det (e_{\a_1},\dots,e_{\a_{n-1}}, e_{\a_{n+1}})=-\sgn(\sigma^{n-i}) \det(e_{\a_{i+1}},\dots,e_{\a_{n-1}},e_{\a_n},e_{\a_1},\dots,e_{\a_{i-1}},e_{\a_{n+1}}),
\end{align*}
where $\a_{n+1} \in \seq{1}{n} \setminus \{\seqalpha{\a}{1}{n-1}\}$
which is true since
\begin{align*}
&\det(e_{\a_{i+1}},\dots,e_{\a_{n-1}},e_{\a_n},e_{\a_1},\dots,e_{\a_{i-1}},e_{\a_{n+1}}) \\
&=\det(e_{\a_{i+1}},\dots,e_{\a_{n-1}},e_{\a_i},e_{\a_1},\dots,e_{\a_{i-1}},e_{\a_{n+1}}) \\
&=(-1)^{(n-1)i-1}\det(e_{\a_{1}},\dots,e_{\a_{n-1}},e_{\a_{n+1}}),
\end{align*}
and
\begin{align*}
\sgn(\sigma^{n-i})=(-1)^{(n-1)(n-i)}.
\end{align*}

2. Lemma \ref{lem:mo3} also gives
\begin{align*}
\left(\ct{m}{n}{F(T^n,\bar e)} \otimes \bar e\right) (e_{\a_1}, \dots, e_{\a_n})
= \frac{1}{n!} \det (e_{\a_1},\dots,e_{\a_{n-1}},(e_1+\dots+e_n)).
\end{align*}
Together with \eqref{eq:1}, we have
\begin{align*}
\left(\sum_{u\in \MN(T^n)}  \ct{m}{n}{F(T^n,u)} \otimes u \right) (e_{\a_1}, \dots, e_{\a_n}) =0,
\end{align*}
if there are $i,j \in \{1,\dots,n\}$ such that $\a_i=\a_j$; and
\begin{align*}
\left(\sum_{u\in \MN(T^n)}  \ct{m}{n}{F(T^n,u)} \otimes u \right)  (e_{\a_1}, \dots, e_{\a_n})
&= \frac{1}{n!} \det(e_{\a_1},\dots,e_{\a_{n-1}},e_{\a_n})\\
&= V(T^n)\lct_{\a_1,\dots,\a_n}
\end{align*}
for pairwise distinct $\a_1,\dots,\a_n$.
Thus \eqref{eq:MR1} holds.
\end{proof}

We close this section by computing valuations in Theorem \ref{thm:contrav} on $sT^n$ for $s>0$.
Observe that $\mathcal {N}_o (sT^n) = \{\bar e\}$. Direct calculations show
\begin{equation}\label{eq615-0}
\begin{aligned}
\ab{\sum_{u \in \mathcal {N}_o (sT^n)} h_{sT^n}^{-p}(u) \zeta\big(V(sT^n,u)\big) \symtp{u}{p}}_{\a_1\dots\a_p}
&= \big(h_{sT^n}^{-p} (\bar e) \zeta(V(sT^n, \bar e)) \symtp{(\bar e)}{p} \big)_{\a_1\dots\a_p}\\
&= s^{-p} \zeta\ab{\frac{s^n}{n!}}.
\end{aligned}
\end{equation}
Also, Lemma \ref{lem:mo3} gives
\begin{equation}\label{eq615-1}
\begin{aligned}
\ct{m}{n+1}{sT^{n}}_{\a_1 \dots \a_{n-1}}
&=\frac{1}{(n+1)!}s^{n+1} \det(e_{\a_1},\dots,e_{\a_{n-1}},e_1+\dots+e_n) \\
&=\frac{1}{(n+1)!}s^{n+1} \lct_{\a_1 \dots  \a_{n-1}} 
\end{aligned}
\end{equation}
and
\begin{equation*}
\begin{aligned}
\ab{\ct{m}{n}{F(sT^n,\bar e)} \otimes \bar e}_{\a_1 \dots  \a_{n}}
=&\frac{1}{n!}s^{n} \det(e_{\a_1},\dots,e_{\a_{n-1}},e_1+\dots+e_n)\\
=&\frac{1}{n!}s^{n} \lct_{\a_1 \dots \a_{n-1}}. 
\end{aligned}
\end{equation*}
Hence
\begin{equation}\label{eq615-2}
\begin{aligned}
\ab{\sigma^r\ab{\ct{m}{n}{F(sT^n,\bar e)} \otimes \bar e}}_{\a_1\dots\a_n}
=&\ab{\ct{m}{n}{F(sT^n,\bar e)} \otimes \bar e}_{\a_{\sigma^{-r}(1)} \dots \a_{\sigma^{-r}(n)}} \\
=&\frac{1}{n!} \lct_{\a_{\sigma^{-r}(1)} \dots \a_{\sigma^{-r}(n-1)}}.
\end{aligned}
\end{equation}

\section{The case $n=p=2$}\label{sec:n=2}
The aim of this section is to prove Theorem \ref{thm:tensor:p} for $n=2$.
It suffices to determine $Z(sT^d)$ for every $s>0$ and $d=0,1,2$ due to Lemma \ref{lemuq} and Theorem \ref{thm:contrav}.
Note that $\rho \phi \rho^{-1} = \phi^{-t}$ for every $\phi \in \SL{2}$.
Thus a mapping $Z:\MP_o^2 \to \tenset{2}{2}$ is an $\SL{2}$ covariant valuation if and only if $Z \circ \rho :\MP_o^2 \to \tenset{2}{2}$ is an $\SL{2}$ contravariant valuation.
Classifications of $\SL{2}$ covariant tensor valuations were established by Wang \cite[Theorem 1.3]{Wang2022} although the valuation $P \mapsto (-1)^{\dim P}V_0(\{o\} \cap \relint P) \lct$ was unfortunately missing (the same omission appears in Theorem 3 in Zeng and Zhou \cite{ZZ2024}). It is a very minor mistake which does not affect the correctness of the rest of the results in Wang \cite{Wang2022},  in particular his Lemma 3.1, Corollary 3.1 and Corollary 3.2, which we rely on (via the rotation $\rho$) to  prove Theorem \ref{thm:tensor:p} for $n=2$.

\begin{lem}\label{thm:wl3p1} \textup{\cite[Lemma 3.1]{Wang2022}}
If $Z: \MPo2 \to \tenset{2}{2}$ is an $\SL{2}$ contravariant valuation, then
there is a constant $d_0\in\R$ such that
\[Z(\set{o})=-d_0\rho
=d_0\left(\begin{array}{cc}
	0 & -1\\
	1 & 0
\end{array}\right).\]
\end{lem}

\begin{lem}\label{thm:wc3p1}\textup{\cite[Corollary 3.1]{Wang2022}}
If $Z: \MPo2 \to \tenset{2}{2}$ is an $\SL{2}$ contravariant valuation, then
there are constants $d_1,d_2\in\R$ such that
\[Z(sT^1)=d_1\dten{2}(s \rho T^1)-d_2\rho =2d_1s^2\left(\begin{array}{cc}
	0 & 0\\
	0 & 1
\end{array}\right)+d_2\left(\begin{array}{cc}
	0 & -1\\
	1 & 0
\end{array}\right)\]
for every $s>0$.
\end{lem}

\begin{lem}\label{thm:wc3p2}\textup{\cite[Corollary 3.2]{Wang2022}}
	If $Z: \MPo2 \to \tenset{2}{2}$ is an $\SL{2}$ contravariant valuation, then
	there are Cauchy functions ${\tilde \zeta},{\tilde \eta}:[0,\infty)\to\R$ and constants $d_1,d_2,d_3,d_4\in\R$ such that
	\begin{align*}
		Z(sT^2)=&d_3M^{2,0}(s \rho T^2)+d_1 \dten{2}(s \rho T^2)+d_4F(s\rho T^2)+G_{\tilde \eta}(s\rho T^2)+H_{\tilde \zeta}(s\rho T^2)-d_2\rho\\
		=&\frac1{24}d_3s^4\left(\begin{array}{cc}
			2 & -1\\
			-1 & 2
		\end{array}\right)+d_1s^2\left(\begin{array}{cc}
			1 & 0\\
			0 & 1
		\end{array}\right)+d_4s^2\left(\begin{array}{cc}
			-1 & 0\\
			0 & 1
		\end{array}\right)\\
		&+{\tilde \eta}(s^2)\left(\begin{array}{cc}
			0 & 1\\
			-1 & 0
		\end{array}\right)+\frac1{s^2}{\tilde \zeta}(s^2)\left(\begin{array}{cc}
			1 & 1\\
			1 & 1
		\end{array}\right)+d_2\left(\begin{array}{cc}
			0 & -1\\
			1 & 0
		\end{array}\right)
	\end{align*}
	for every $s>0$.
\end{lem}
The mappings $F$, $G_{\tilde \eta}$ and $H_{\tilde \zeta}$ are matrix valuations on $\MPo2$ introduced in \cite{Wang2022} defined by the following.
\[F(P):=\left\{\begin{array}{ll}
	v^2-w^2, & \text{if }\dim P=2\text{ and }P\text{ has two edges }[o,v]\text{ and }[o,w]\text{ with }\det(v,w)>0;\\
	0, & \text{otherwise};
\end{array}\right.\]
\[H_{\tilde \zeta}(P):=\sum_{i=2}^r  \frac{{\tilde \zeta}(\det(v_{i-1},v_i))}{\det(v_{i-1},v_i)^2}  (v_{i-1}-v_i)^2\]
if $\dim P=2$ and $P=[o,v_1,\dots,v_r]$ with $o\in\bd P$ and the vertices $o,v_1,\dots,v_r$ are labeled counter-clockwise;
\[H_{\tilde \zeta}(P):=\frac{{\tilde \zeta}(\det(v_r,v_1))}{\det(v_r,v_1)^2}  (v_r-v_1)^2
+\sum_{i=2}^r  \frac{{\tilde \zeta}(\det(v_{i-1},v_i))}{\det(v_{i-1},v_i)^2}  (v_{i-1}-v_i)^2\]
if $o\in\intr P$ and $P=[v_1,\dots,v_r]$ with the vertices $v_1,\dots,v_r$ labeled counter-clockwise;
\[H_{\tilde \zeta}(P):=0\]
if $P=\set{o}$ or $P$ is a line segment; and
\[G_{\tilde \eta}(P):=\sum_{i=2}^r  \frac{{\tilde \eta}(\det(v_{i-1},v_i))}{\det(v_{i-1},v_i)}  (v_{i-1}\otimes v_i-v_i\otimes v_{i-1})\]
if $\dim P=2$ and $P=[o,v_1,\dots,v_r]$ with $o\in\bd P$ and the vertices $o,v_1,\dots,v_r$ are labeled counter-clockwise;
\[G_{\tilde \eta}(P):=\frac{{\tilde \eta}(\det(v_r,v_1))}{\det(v_r,v_1)}  (v_r\otimes v_1-v_1\otimes v_r)
+\sum_{i=2}^r  \frac{{\tilde \eta}(\det(v_{i-1},v_i))}{\det(v_{i-1},v_i)}  (v_{i-1}\otimes v_i-v_i\otimes v_{i-1})\]
if $o\in\intr P$ and $P=[v_1,\dots,v_r]$ with the vertices $v_1,\dots,v_r$ labeled counter-clockwise;
otherwise,
\[G_{\tilde \eta}(P):=0.\]

Now we prove Theorem \ref{thm:tensor:p} for $n=2$. We restate it here.
\begin{thm}\label{thm:p=n=2}
A mapping $Z: \MPo2 \to \tenset{2}{2}$ is an $\SLn$ contravariant valuation if and only if there are Cauchy functions $\zeta,\eta : [0,\infty) \to \R$ and constants $a,b$, $c_0',c_0, c_1$ such that
\begin{equation}\label{eq904-1}
\begin{split}
ZP&= \sum_{u \in \mathcal {N}_o (P)} h_P^{-2}(u) \zeta(V(P,u)) \symtp{u}{2}
+ c_1  \sum_{u \in \MN_o(P)} \ct{m}{2}{F(P,u)} \otimes u\\
&\qquad + \big(c_0'(-1)^{\dim P}V_0(\{o\} \cap \relint P) + c_0V_0(P) + \eta(V_2(P))\big) \lct \\
&\qquad + a  M^{2,0}(\rho P)
+ b \dten{2} (\rho P)
\end{split}
\end{equation}
for every $P \in \MP_o^2$.
\end{thm}

Remark: Since the previous valuations $P \mapsto d_4F(P)+G_{\tilde \eta}(P)+H_{\tilde \zeta}(sP)-d_2\rho$ are apparently restricted to planar domains,
the challenge here actually lies in conjecturing their appropriate modifications to be suitable also for higher dimensions.
The subsequent verification in the following is straightforward.

\begin{proof}[Proof of Theorem \ref{thm:p=n=2}.]
Set
	\[c_0'=-d_0,\quad b=d_1,\quad c_0=-d_2,\quad a=d_3,\quad c_1=-2d_4,\]
	\[\zeta=2{\tilde \zeta},\qquad \eta(t)=2{\tilde \eta}(t)+2d_4 t,\quad t\geq0\]
in Lemmas \ref{thm:wl3p1}, \ref{thm:wc3p1} and \ref{thm:wc3p2}.
Direct calculations show that
\[Z(\{o\})=c_0' \ab{\begin{array}{cc} 0 & 1 \\ -1 & 0\end{array}},\qquad
Z(sT^1)=c_0\ab{\begin{array}{cc} 0 & 1 \\ -1 & 0\end{array}} + 2bs^2 \ab{\begin{array}{cc} 0 & 0 \\ 0 & 1\end{array}},\]
and
\begin{align*}
	Z(sT^2)
	=&\frac{1}{s^2}\zeta\ab{\frac{s^2}{2}} \ab{\begin{array}{cc} 1 & 1 \\ 1 & 1 \end{array}}
	+c_1 \frac{s^2}{2}\ab{\begin{array}{cc}1 & 1 \\ -1 & -1\end{array}}
	+\ab{c_0+\eta\ab{\frac{s^2}{2}}} \ab{\begin{array}{cc} 0 & 1 \\ -1 & 0\end{array}}\\
	&+\frac{1}{24}a s^4\ab{\begin{array}{cc}2 &-1 \\ -1 &2\end{array}}
	+bs^2 \ab{\begin{array}{cc} 1 &0 \\ 0 & 1\end{array}}
\end{align*}
for every $s>0$.
Therefore,
\begin{align}
	ZP&= \sum_{u \in \mathcal {N}_o (P)} h_P^{-2}(u) \zeta(V(P,u)) \symtp{u}{2}
	+ c_1  \sum_{u \in \MN_o(P)} \ct{m}{2}{F(P,u)} \otimes u\\
	&\qquad + \big(c_0'(-1)^{\dim P}V_0(\{o\} \cap \relint P) + c_0V_0(P) + \eta(V_2(P))\big) \lct \\
	&\qquad + a  M^{2,0}(\rho P) 
	+ b \dten{2} (\rho P)
\end{align}
for $P=\{o\}$, and for $P=sT^1$ and $P=sT^2$ with any $s>0$.
Hence the desired result follows from Lemma \ref{lemuq} and Theorem \ref{thm:contrav}.
\end{proof}

\section{The case \texorpdfstring{$n>2$}{n>2}}\label{sec:n>2}
In this section, we always assume that $n \ge 3$ and $p \ge 2$ are such that $p \le n$.
The aim of this section is to prove the following three theorems.
\begin{thm}\label{thm:low}
Let $n \ge 3$ and $p \ge 2$ such that $p \le n$.
If $Z \in \tenval(\MPon;\tenset{n}{p})$,
then there exist constants $c_0',c_0 \in \R$ such that
\begin{align*}
Z(sT^d)=
\begin{cases}
\big(c_0'(-1)^{d}V_0(\{o\} \cap \relint sT^d) + c_0V_0(sT^d) \big) \lct, & \text{if~} p=n,\\
0, & \text{if~} p<n
\end{cases}
\end{align*}
for every $0 \le d \le n-1$ and $s>0$.
\end{thm}

\begin{thm}\label{thm:tensor:po}
Let $n \ge 3$ and $p \ge 2$ such that $p \le n$.
If $Z^+ \in \tenval_0(\MPon;\tenset{n}{p})$ is simple, then there exists a Cauchy function $\zeta: [0,\infty) \to \R$ such that
\begin{align*}
Z^+(sT^n)= \sum_{u \in \mathcal {N}_o (sT^n)} h_{sT^n}^{-p}(u) \zeta\big(V(sT^n,u)\big) \symtp{u}{p}
\end{align*}
for every $s>0$.
\end{thm}

\begin{thm}\label{thm:tensor:neg}
Let $n \ge 3$ and $p \ge 2$ such that $p \le n$.
If $Z^- \in \tenval_1(\MPon;\tenset{n}{p})$ is simple, then there exist a Cauchy function $\eta: [0,\infty) \to \R$ and constants $c_1,\dots,c_{n-1} \in\R$ such that
\begin{align*}
Z^- (sT^n)=
\begin{cases}
c \ct{m}{n+1}{sT^n}, & \text{if } p=n-1,\\
\sum_{r=0}^{n-2} c_{r+1} \sigma^r \left(\sum_{u \in \MN_o(sT^n)} \ct{m}{n}{F(sT^n,u)} \otimes u \right) +\eta\big(V_n(sT^n)\big) \lct, & \text{if } p=n
\end{cases}
\end{align*}
for every $s>0$.
\end{thm}

Note that Theorem \ref{thm:contrav} proves the ``if" parts of Theorems \ref{thm:tensor:p+2}, \ref{thm:tensor:p+1} and \ref{thm:tensor:p}.
Furthermore, \S \ref{sec:DV}, Lemma \ref{lemuq} and Theorem \ref{thm:contrav} show that Theorems \ref{thm:low}, \ref{thm:tensor:po} and \ref{thm:tensor:neg} are sufficient to prove the ``only if" parts of Theorems \ref{thm:tensor:p+2}, \ref{thm:tensor:p+1} and \ref{thm:tensor:p} for $n > 2$.
In fact, supposing $Z \in \tenval(\MPon;\tenset{n}{p})$, Lemma \ref{lemuq} for the lower-dimensional case, Theorems \ref{thm:contrav} and \ref{thm:low} show that
\begin{align*}
ZP=
\begin{cases}
\big(c_0'(-1)^{d}V_0(\{o\} \cap \relint P) + c_0V_0(P) \big) \lct, & \text{~if~} p=n,\\
0, & \text{~if~} p<n
\end{cases}
\end{align*}
for all $P \in \MPon$ with $\dim P<n$.
Consider $\tilde ZP:=ZP-\big(c_0'(-1)^{\dim P}V_0(\{o\} \cap \relint P) + c_0V_0(P)\big) \lct$ if $p=n$ and $\tilde ZP:=ZP$ if $p<n$ for all $P \in \MPon$.
Thus $\tilde Z$ is a simple $\sln$ contravariant valuation on $\MPon$,
and then \S \ref{sec:DV}, Lemma \ref{lemuq}, Theorems \ref{thm:contrav}, \ref{thm:tensor:po} and \ref{thm:tensor:neg} determine $\tilde Z$.

\vskip 10pt
\textbf{Induction step.}
The proofs of Theorems \ref{thm:low}, \ref{thm:tensor:po} and \ref{thm:tensor:neg} use the following induction step.
Let $Z:\MPon \to \tenset{n}{p}$ be a mapping.
For real $s>0$, $q \in \{1,\dots, p\}$, $1 \leq k_1<\dots<k_q  \le p$ and $i \in \{1,\dots,n\}$, we define a mapping $\indval{s}{i}{k_1,\dots,k_q}$ that maps $\MP_o(e_i^{\bot})=\{P \in \MPon: P \subset e_i^{\bot}\}$ to $(e_i^{\bot})^{\otimes_{(p-q)}}$ (recall $(e_i^{\bot})^{\otimes_0}=\R$) by
\begin{align}
\indval{s}{i}{k_1,\dots,k_q}(P)_{\a_1 \dots \a_{\hat{k}_1} \dots \a_{\hat{k}_q}\dots \a_p}:=Z([P,se_i])_{\a_1 \dots \a_{k_1} \dots \a_{k_q} \dots \a_p}, ~P\in \MP_o(e_i^{\bot})
\end{align}
for any $\a_1 ,\dots,  \a_{\hat{k}_1} ,\dots, \a_{\hat{k}_q} ,\dots, \a_p \in \{1,\dots \hat{i} ,\dots, n\}$ and $\a_{k_1} =\a_{k_2}=\dots =\a_{k_q}=i$.
Here $\a_{\hat{k}_1}, \dots, \a_{\hat{k}_q}$ means that subindices $\hat{k}_1, \dots, \hat{k}_q$ do not appear (only $p-q$ indices are left).
If $Z$ is an $\slpm{n}$-$\delta$-contravariant valuation, then $\indval{s}{i}{k_1,\dots,k_q}$ is an $\slpm{e_i^{\bot}}$-$\delta$-contravariant valuation.
Here $\slpm{e_i^{\bot}} = \{\phi \in \SL{e_i^{\bot}} : \det_{n-1} \phi = 1 \text{~or~} -1\}$ (Recall that $e_i^{\bot} \cong \R^{n-1}$).
Briefly for $q=1$, we write $\indval{s}{i}{k}$ instead of $\indval{s}{i}{k_1}$, i.e.,
\begin{align}\label{ind1}
\indval{s}{i}{k}(P)_{\a_1 \dots \a_{\hat{k}} \dots \a_p}=Z([P,se_i])_{\a_1 \dots \a_{k} \dots \a_p},
\end{align}
for any $\a_1, \dots,  \a_{\hat{k}}, \dots, \a_p \in \{1,\dots \hat{i} ,\dots, n\}$ and $\a_{k} =i$.
Here $\a_{\hat{k}}$ means that the subindex $\hat{k}$ does not appear.

Let $Z \in \tenval(\MPon;\tenset{n}{p})$, $q \in \seq{1}{p}$ and $\a_{k_1} =\a_{k_2}=\dots =\a_{k_q}=i$ while $\a_1, \dots,  \a_{\hat{k}_1} \dots ,\a_{\hat{k}_q}, \dots ,\a_p \in \{1,\dots \hat{i} ,\dots, n\}$.
Set $\phi \in \sln$ with $\phi e_i = se_i$ and $\phi e_j=s^{-1/(n-1)}$ for all $j \in \{1,\dots \hat{i} ,\dots, n\}$.
Then
\begin{equation}\label{eq:hom}
\begin{aligned}
Z(sT^d)_{\a_1 \dots \a_{k_1} \dots \a_{k_q} \dots \a_p}
&=Z(\phi[s^{n/(n-1)}\hat T^{d-1}_i,e_i])_{\a_1 \dots \a_{k_1} \dots \a_{k_q} \dots \a_p} \\
&=s^{-q+(p-q)/(n-1)} Z([s^{n/(n-1)}\hat T^{d-1}_i,e_i])_{\a_1 \dots \a_{k_1} \dots \a_{k_q} \dots \a_p}
\end{aligned}
\end{equation}
for every $1 \le d \le n$ with $i \le d$.
Therefore,
\begin{equation}\label{eq:ind1a}
\begin{aligned}
Z(sT^d)_{\a_1 \dots \a_{k_1} \dots \a_{k_q} \dots \a_p}
=s^{-q+(p-q)/(n-1)} \indval{1}{i}{k_1,\dots,k_q}(s^{n/(n-1)}\hat{T}_i^{d-1})_{\a_1 \dots \a_{\hat{k}_1} \dots \a_{\hat{k}_q}\dots \a_p}.
\end{aligned}
\end{equation}

\textbf{Induction hypothesis.} Let $n>3$. Assume that the ``only if" parts of Theorems \ref{thm:tensor:p+2}, \ref{thm:tensor:p+1} and \ref{thm:tensor:p} hold on $\R^{n-1}$ for all tensor valuations of order $q$ with $2 \le q < p$.

For the lower-dimensional simplices, the induction hypothesis is equivalent to assuming that Theorem \ref{thm:low} holds on $\R^{n-1}$ for all tensor valuations of order $q$ with $2 \le q < p$.
Together with Theorem \ref{thm:p=0}, Theorem \ref{thm:p=1} and \eqref{eq:ind1a}, we get
\begin{equation}\label{eq:indlow1}
\begin{aligned}
&Z(sT^d)_{\a_1 \dots \a_{k_1} \dots \a_{k_q} \dots \a_p}=s^{-q+(p-q)/(n-1)} \indval{1}{i}{k_1,\dots,k_q}(s^{n/(n-1)}\hat{T}_i^{d-1})_{\a_1 \dots \a_{\hat{k}_1} \dots \a_{\hat{k}_q}\dots \a_p}\\
&=\begin{cases}
c_0(\indfconst{i}{k_1,\dots,k_q}) s^{-q+1} \lct^{e_i^\bot}_{\a_1 \dots \a_{\hat{k}_1} \dots \a_{\hat{k}_q}\dots \a_p}, & \text{~if~} p-q=n-1,\\
c_0(\indfconst{i}{k_1,\dots,k_q}) s^{-q}, & \text{~if~} p=q,\\
0, & \text{~otherwise}
\end{cases}
\end{aligned}
\end{equation}
for $2 \le d \le n-1$ and $1 \le i \le d$;
\begin{equation}\label{eq:indlow2}
\begin{aligned}
&Z(sT^1)_{\a_1 \dots \a_{k_1} \dots \a_{k_q} \dots \a_p}=s^{-q+(p-q)/(n-1)} \indval{1}{i}{k_1,\dots,k_q}(s^{n/(n-1)}\{o\})_{\a_1 \dots \a_{\hat{k}_1} \dots \a_{\hat{k}_q}\dots \a_p}\\
&=\begin{cases}
\big(c_0' (\indfconst{i}{k_1,\dots,k_q}) + c_0(\indfconst{i}{k_1,\dots,k_q}))s^{-q+1} \lct^{e_i^\bot}_{\a_1 \dots \a_{\hat{k}_1} \dots \a_{\hat{k}_q}\dots \a_p}, & \text{~if~} p-q=n-1,\\
\big(c_0' (\indfconst{i}{k_1,\dots,k_q}) + c_0(\indfconst{i}{k_1,\dots,k_q}))s^{-q}, & \text{~if~} p=q,\\
0, & \text{~otherwise},
\end{cases}
\end{aligned}
\end{equation}
for $i=1$, where $c_0(\indfconst{i}{k_1,\dots,k_q}), ~ c_0'(\indfconst{i}{k_1,\dots,k_q}) \in \R$ which depend on $i$ and $k_1,\dots,k_q$, but not on $d$.
For the case $p-q=n-1$, we have $q=1$. Therefore, we will write $c_0(\indfconst{i}{k_1})$ and $c_0'(\indfconst{i}{k_1})$ as $c_0(\indfconst{i}{k})$ and $c_0'(\indfconst{i}{k})$.
For the case $p=q$, we have $k_l=l$ for every $1 \le l\le p$.

For the full-dimensional simplices, decompose $Z = Z^+ + Z^-$ with $Z^+ \in \tenval_0(\MPon;\tenset{n}{p})$ and $Z^- \in \tenval_1(\MPon;\tenset{n}{p})$ as in \S \ref{sec:DV}.
Applying the induction hypothesis to $Z^+ \in \tenval_0(\MPon;\tenset{n}{p})$ and $Z^- \in \tenval_1(\MPon;\tenset{n}{p})$, and using Theorems \ref{thm:p=0} and \ref{thm:p=1} together with \eqref{eq615-0}, \eqref{eq615-1}, \eqref{eq615-2}, \eqref{eq:ind1a} and Theorem \ref{thm:contrav}, we can find Cauchy functions $\zeta(\indf{i}{k_1,\dots,k_q}{\cdot})$ and $\eta(\indf{i}{k}{\cdot})$ on $[0,\infty)$ and $c(\indfconst{i}{\seqnb{k_1}{k_q}}),c_{r}(\indfconst{i}{k}) \in \R$ for $r=0,\dots,n-2$ depending on $i$ and $k_1,\dots,k_q$ such that
\begin{equation}\label{eq:ind2}
\begin{aligned}
Z^+(sT^n)_{\a_1 \dots \a_{k_1} \dots \a_{k_q} \dots \a_p}
&=s^{-q+(p-q)/(n-1)}\indval{1}{i}{\seqnb{k_1}{k_q}}^+(s^{n/(n-1)}\hat T^{n-1}_i)_{\a_1 \dots \a_{\hat{k}_1} \dots \a_{\hat{k}_q}\dots \a_p}  \\
&=\zeta(\indf{i}{k_1,\dots,k_q}{s^n}) s^{-(p-q)n/(n-1)-q+(p-q)/(n-1)} \\
&=\zeta(\indf{i}{k_1,\dots,k_q}{s^n}) s^{-p}
\end{aligned}
\end{equation}
for all $1 \le q < p \le n$;
\begin{equation}\label{eq:ind2a}
\begin{aligned}
Z^+(sT^n)_{\a_1 \dots \a_{k_1} \dots \a_{k_q} \dots \a_p}
&=\ab{\zeta(\indf{i}{1,\dots,p}{s^n}) + c_0(\indfconst{i}{1,\dots,p})}s^{-n}
\end{aligned}
\end{equation}
for $q=p$;
\begin{equation}\label{eq:ind2a-}
\begin{aligned}
Z^-(sT^n)_{\a_1 \dots \a_{k_1} \dots \a_{k_q} \dots \a_p}
&=s^{-q+(p-q)/(n-1)}\indval{1}{i}{\seqnb{k_1}{k_q}}^-(s^{n/(n-1)}\hat T^{n-1}_i)_{\a_1 \dots \a_{\hat{k}_1} \dots \a_{\hat{k}_q}\dots \a_p}\\
&= c(\indfconst{i}{\seqnb{k_1}{k_q}})s^{(n^2+p-qn)/(n-1)} \lct_{\a_1 \dots \a_{\hat{k}_1} \dots \a_{\hat{k}_q}\dots \a_p}^{e_i^\bot}
\end{aligned}
\end{equation}
for $p=n-1$ with $q=1$ or for $p=n$ with $q=2$;
\begin{equation}\label{eq:ind2b-}
\begin{aligned}
Z^-(sT^n)_{\a_1 \dots \a_{k} \dots \a_n}
=\sum_{r=0}^{n-3} c_{r+1}(\indfconst{i}{k})s^{n}
\lct_{\a_{\sigma^{-r}(1)} \dots \a_{\hat{k}} \dots \a_{\sigma^{-r}(n-1)}}^{e_i^\bot}
+ \big(\eta(\indf{i}{k}{s^n}) + c_0(\indfconst{i}{k}) \big) \lct_{\a_1 \dots \a_{\hat{k}} \dots \a_n}^{e_i^\bot}
\end{aligned}
\end{equation}
for $p=n$ with $q=1$ and $k \neq n$, where $\sigma=\ab{\begin{smallmatrix}	1  & \cdots &k-1 &k+1 &\cdots  & n-1 & n \\ 2  & \cdots & k+1 & k+2 &\cdots & n & 1 \end{smallmatrix}}$;
\begin{equation}\label{eq:ind2bbb-}
\begin{aligned}
Z^-(sT^n)_{\a_1 \dots \a_n}
=&\sum_{r=0}^{n-3} c_{r+1}(\indfconst{i}{n})s^{n}
\lct_{\a_{\sigma^{-r}(1)} \dots \a_{\sigma^{-r}(n-2)}}^{e_i^\bot}
+ \big(\eta(\indf{i}{n}{s^n}) + c_0(\indfconst{i}{n}) \big) \lct_{\a_1 \dots \a_{n-1}}^{e_i^\bot}
\end{aligned}
\end{equation}
for $p=n$ with $q=1$ and $k = n$,
where $\sigma=\ab{\begin{smallmatrix}	1  & \cdots  & n-2 & n-1 \\ 2  & \cdots   & n-1 & 1 \end{smallmatrix}}$; and
\begin{equation}\label{eq:ind2c-}
\begin{aligned}
Z^-(sT^n)_{\a_1 \dots \a_{k_1} \dots \a_{k_q} \dots \a_p}=0,
\end{aligned}
\end{equation}
otherwise.
Note that $c_0(\indfconst{i}{k})$ in \eqref{eq:indlow1} and \eqref{eq:indlow2} for the case $p-q=n-1$ is the same as $c_0(\indfconst{i}{k})$ in \eqref{eq:ind2b-} and \eqref{eq:ind2bbb-},
and $c_0(\indfconst{i}{1,\dots,p})$ in \eqref{eq:indlow1} and \eqref{eq:indlow2} for the case $p=q$ and $c_0(\indfconst{i}{1,\dots,p})$ in \eqref{eq:ind2a} are the same; since they come from applying induction to the same valuation $P \mapsto V_0(P) \lct$.

We can simplify \eqref{eq:ind2b-} and \eqref{eq:ind2bbb-} if we further assume $\a_1, \dots, \a_n$ are pairwise distinct numbers in $\{1,\dots, n\}$ with $\a_{k}=i$.
Assume first $k \neq n$. Observe that
\begin{align*}
&\lct_{\a_{\sigma^{-r}(1)} \dots \a_{\sigma^{-r}(k-1)} \a_{\sigma^{-r}(k+1)} \dots \a_{\sigma^{-r}(n-1)}}^{e_i^\bot}\\
=&\det\nolimits_{(n-1)} \ab{e_{\a_{\sigma^{-r}(1)}},\dots,e_{\a_{\sigma^{-r}(k-1)}},e_{\a_{\sigma^{-r}(k+1)}},\dots,e_{\a_{\sigma^{-r}(n)}}}\\
=&(-1)^{nr} \det\nolimits_{(n-1)} \ab{e_{\a_1},\dots,e_{\a_{k-1}},e_{\a_{k+1}},\dots,e_{\a_n}}
\end{align*}
and
\begin{align*}
\lct_{\a_1 \dots \a_{\hat{k}} \dots \a_n}^{e_i^\bot}
&= \det\nolimits_{(n-1)} \ab{e_{\a_1},\dots,e_{\a_{k-1}},e_{\a_{k+1}},\dots,e_{\a_n}} \\
&=(-1)^{i+k} \lct_{\a_1 \dots \a_{k} \dots \a_n}.
\end{align*}
By \eqref{eq:ind2b-}, we have
\begin{equation}\label{eq:ind2d-}
\begin{aligned}
&Z^-(sT^n)_{\a_1 \dots \a_{k} \dots \a_n}\\
=&(-1)^{i+k} \ab{\sum_{r=0}^{n-3} (-1)^{nr} c_{r+1}(\indfconst{i}{k})s^n +\eta(\indf{i}{k}{s^n}) +  c_0(\indfconst{i}{k})} \lct_{\a_1 \dots \a_{k} \dots \a_n},
\end{aligned}
\end{equation}
where $\a_1, \dots, \a_n$ are pairwise distinct numbers in $\{1,\dots, n\}$ with $\a_{k}=i$.
Similarly we obtain from \eqref{eq:ind2bbb-} that \eqref{eq:ind2d-} also holds for $k=n$.

For $n=3$, the induction relies on the planar case (Theorems \ref{thm:p=0}, \ref{thm:p=1} and \ref{thm:p=n=2}), where additional valuations appear:
\begin{align*}
M^{2,0}(\rho sT^2)=\frac{1}{24} s^4(2e_1^2-(e_1\otimes e_2+e_2\otimes e_1)+2e_2^2), ~M^{2,0}(\rho sT^1)=0,
\end{align*}
\begin{align*}
\dten{2} (\rho sT^2)=s^2 (e_1^2+e_2^2), ~\dten{2} (\rho sT^1)= 2 s^2 e_2,
\end{align*}
and
\begin{align*}
\dten{1}(\rho sT^2)=s(e_1-e_2), ~\dten{1}(\rho sT^1)=-2s e_2.
\end{align*}
But we will show in Lemma \ref{lem:n=3} that these valuations have no corresponding valuations in $3$-dimensional space.

\begin{lem}\label{lem:n=3}
Let $p=2,3$ and $Z \in \tenval(\MP_o^3;\tenset{3}{p})$.
Decompose $Z= Z^+ + Z^-$ for $Z^+ \in \tenval_0(\MP_o^3;\tenset{3}{p})$ and $Z^- \in \tenval_1(\MP_o^3;\tenset{3}{p})$ as in \S \ref{sec:DV}.
Then \eqref{eq:indlow1}, \eqref{eq:indlow2}, \eqref{eq:ind2}, \eqref{eq:ind2a}, \eqref{eq:ind2a-}, \eqref{eq:ind2b-}, \eqref{eq:ind2bbb-}, \eqref{eq:ind2c-} and \eqref{eq:ind2d-} still hold for $n=3$.
\end{lem}
\begin{proof}
Let $s>0$.
First, if $\a_1=\dots=\a_p=i$, then \eqref{eq:ind1a} for both $Z^+ \in \tenval_0(\MP_o^3;\tenset{3}{p})$ and $Z^- \in \tenval_1(\MP_o^3;\tenset{3}{p})$ together with Theorem \ref{thm:p=0} gives (similar to \eqref{eq:ind2a})
\begin{equation*}
Z^+(sT^3)_{\a_1\dots\a_p}=\big(\zeta(\indf{i}{1,\dots,p}{s^3}) + c_0(\indfconst{i}{1,\dots,p}) \big)s^{-p}
\end{equation*}
and
\begin{equation}\label{eq604-2a}
Z^-(sT^3)_{\a_1\dots\a_p}=0
\end{equation}
for some Cauchy functions $\zeta(\indf{i}{1,\dots,p}{\cdot}) : [0,\infty) \to \R$ and constants $c_0(\indfconst{i}{1,\dots,p})\in \R$.
Thus \eqref{eq:ind2a} and \eqref{eq:ind2c-} hold for this setting.
Similarly, applying induction to $sT^2$ and $sT^1$, we get \eqref{eq:indlow1} and \eqref{eq:indlow2} for this setting.

Second, assume $p=3$.
For fixed $i, j \in \set{1,2,3}$ with $i \neq j$, assume $\a_k=i$ and $\a_{l_1}=\a_{l_2}=j$ for pairwise distinct $k,l_1,l_2\in \set{1,2,3}$. Assume w.l.o.g. $l_1<l_2$.
Then \eqref{eq:ind1a} for $Z^+ \in \tenval_0(\MP_o^3;\tenset{3}{3})$ and $Z^- \in \tenval_1(\MP_o^3;\tenset{3}{3})$, Theorems \ref{thm:contrav} and \ref{thm:p=n=2} give
\begin{align*}
Z^+(sT^3)_{\a_1 \a_2 \a_3}&=\indval{1}{i}{k}^+(s^{3/2}\hat T^{2}_i)_{jj}  \\
&=\zeta(\indf{i}{k}{s^3})s^{-3} + a(\indfconst{i}{k})s^6 M^{2,0}(\rho \hat T^{2}_i)_{jj} + b(\indfconst{i}{k})s^3 \dten{2}(\rho \hat T^{2}_i)_{jj},
\end{align*}
\begin{align*}
Z^-(sT^3)_{\a_1 \a_2 \a_3}
=\indval{1}{i}{k}^-(s^{3/2}\hat T^{2}_i)_{\a_{l_1}\a_{l_2}}
=c_1(\indfconst{i}{k})s^3\lct_{\a_{l_1}}^{e_i^\bot};
\end{align*}
and \eqref{eq:ind1a} for $Z^+ \in \tenval_0(\MP_o^3;\tenset{3}{3})$ and $Z^- \in \tenval_1(\MP_o^3;\tenset{3}{3})$, Theorems \ref{thm:p=1} and \ref{thm:contrav} give
\begin{align*}
Z^+\ab{sT^3}_{\a_1 \a_2 \a_3}
&=s^{-3/2}\indval{1}{j}{l_1,l_2}^+(s^{3/2}\hat T^{2}_j)_{i} =\zeta(\indf{j}{l_1,l_2}{s^3})s^{-3},
\end{align*}
\begin{align*}
Z^-(sT^3)_{\a_1 \a_2 \a_3}
&=s^{-3/2}\indval{1}{j}{l_1,l_2}^-(s^{3/2}\hat T^{2}_j)_{i} \\
&=c(\indfconst{j}{l_1,l_2})s^3 \lct_{i}^{e_j^\bot} + b(\indfconst{j}{l_1,l_2}) \dten{1}(\rho \hat T^{2}_j)_{i}
\end{align*}
for some Cauchy functions $\zeta(\indf{i}{k}{\cdot}), \zeta(\indf{j}{l_1,l_2}{\cdot}):[0,\infty) \to \R$ and constants $a(\indfconst{i}{k}),~ b(\indfconst{i}{k})$, $c_1(\indfconst{i}{k}),~ c(\indfconst{j}{l_1,l_2}),~ b(\indfconst{j}{l_1,l_2}) \in \R$ for all $i,j \in \{1,2,3\}$ and $k,l_1,l_2\in \{1,2,3\}$.
Note that a Cauchy function is rational homogeneous.
Comparing above equations for $sT^3$, we find
\begin{align}\label{eq622-1}
&b(\indfconst{j}{l_1,l_2})=a(\indfconst{i}{k})=b(\indfconst{i}{k})=0,
\end{align}
for all $i,j \in \{1,2,3\}$ and $k,l_1,l_2\in \{1,2,3\}$.
and then \eqref{eq:ind2}, \eqref{eq:ind2a-}, \eqref{eq:ind2b-} and \eqref{eq:ind2bbb-} hold for this setting.
For pairwise distinct $\a_1,\a_2,\a_3$, we can do similar induction to find that  \eqref{eq:ind2}, \eqref{eq:ind2b-} and \eqref{eq:ind2bbb-} also hold by \eqref{eq622-1}.

Let $k,l_1,l_2\in \set{1,2,3}$ be pairwise distinct. Applying induction to $sT^2$ with  \eqref{eq:ind1a}, Theorems \ref{thm:p=1}, \ref{thm:p=n=2} and \ref{thm:contrav}, we have
\begin{align}
Z^+(sT^{2})_{\a_1 \a_2 \a_3}
=\indval{1}{i}{k}^+(s^{3/2}\hat T^{1}_i)_{\a_{l_1} \a_{l_2}} = b(\indfconst{i}{k})\dten{2}([o,-s^{3/2}e_3])_{\a_{l_1} \a_{l_2}},
\end{align}
\begin{align}
Z^-(sT^{2})_{\a_1 \a_2 \a_3}=0,
\end{align}
if $1 \le i \le 2$, $\a_k=i$ and $\a_{l_1},\a_{l_2} \in \{1,2,3\} \setminus \{i\}$; and
\begin{align}
Z^+(sT^{2})_{\a_1 \a_2 \a_3}=0,
\end{align}
\begin{align}
Z^-(sT^{2})_{\a_1 \a_2 \a_3}=\indval{1}{j}{l_1,l_2}^-(s^{3/2}\hat T^{1}_j)_{i} = b(\indfconst{j}{l_1,l_2})\dten{1}([o,-s^{3/2}e_3])_i,
\end{align}
if $1 \le j \le 2$, $\a_{l_1}=\a_{l_2}=j$ and $\a_k \neq j$.
By \eqref{eq622-1}, we conclude \eqref{eq:indlow1} for this setting.
Similarly, applying induction to $sT^1$ concludes \eqref{eq:indlow2} for this setting.

Finally, assume $p=2$. By the first step, we only need to consider the case $\a_1\neq \a_2$.
Assume $\a_k=i$, $\a_l=j$ with $k \neq l$ and $i \neq j$.
Equations \eqref{eq:ind1a} for $Z^+ \in \tenval_0(\MP_o^3;\tenset{3}{2})$ and $Z^- \in \tenval_1(\MP_o^3;\tenset{3}{2})$ together with Theorems \ref{thm:p=1} and \ref{thm:contrav} imply
\begin{align*}
Z^+(sT^3)_{\a_1\a_2}
&=s^{-1/2}\indval{1}{i}{k}^+(s^{3/2}\hat T^{2}_i)_{j}=\zeta(\indf{i}{k}{s^3})s^{-2},
\end{align*}
\begin{align*}
Z^-(sT^3)_{\a_1\a_2}
&=s^{-1/2}\indval{1}{i}{k}^-(s^{3/2}\hat T^{2}_i)_{j} \\
&=c(\indfconst{i}{k})s^4 \lct_{j}^{e_i^\bot}  + b(\indfconst{i}{k})s \dten{1}(\rho \hat T^{2}_i)_{j},
\end{align*}
for some Cauchy function $\zeta(\indf{i}{k}{\cdot}): [0,\infty) \to \R$ and constants $c(\indfconst{i}{k}),b(\indfconst{i}{k}) \in \R$ for all $i \in \{1,2,3\}$ and $k\in \{1,2\}$.
Thus \eqref{eq:ind2} holds for this setting.

Further by $Z^- \in \tenval_1(\MP_o^3;\tenset{3}{2})$, we have
\begin{align*}
Z^-(sT^3)_{12}=-Z^-(sT^3)_{21}=Z^-(sT^3)_{31}=-Z^-(sT^3)_{13}.
\end{align*}
Note that
\begin{align*}
\dten{1}(\rho \hat T^{2}_1)=e_2-e_3, ~\dten{1}(\rho \hat T^{2}_2)=e_1-e_3, ~\dten{1}(\rho \hat T^{2}_3)=e_1-e_2.
\end{align*}
Hence
\begin{align*}
&c(1;1)s^4 + b(1;1)s=-c(1;2)s^4 - b(1;2)s\\
=&-c(2;1)s^4 - b(2;1)s=c(2;2)s^4 + b(2;2)s\\
=&c(3;1)s^4 + b(3;1)s=-c(3;2)s^4 - b(3;2)s.
\end{align*}
In conclusion,
\begin{align*}
c(i;k)=(-1)^{i+k} c,  ~~b(i;k)=(-1)^{i+k} b
\end{align*}
with some constants $b=b(1;1), c=c(1;1) \in \R$ and
\begin{align}\label{eq604-4a}
Z^-(sT^3)_{ij}=
(-1)^{i+1} (c s^4 + b s)\lct_{j}^{e_i^\bot}.
\end{align}

We need to show that $b=0$.
Applying induction to $sT^2$ with  \eqref{eq:ind1a}, Theorems \ref{thm:p=1}, \ref{thm:p=n=2} and \ref{thm:contrav}, we have
\begin{align}
Z^+(sT^{2})_{\a_1 \a_2}=0,
\end{align}
\begin{align}\label{eq604-4aaa}
Z^-(sT^{2})_{\a_1 \a_2}=\indval{1}{i}{k}^-(s^{3/2}\hat T^{1}_i)_{j}=b(\indfconst{i}{k})\dten{1}([o,-s^{3/2}e_3])_j
\end{align}
if $\a_k=i$ and $\a_l=j$ with $k \neq l$, $i \neq j$ and $1 \le i \le 2$.
Since $Z^- \in \tenval_1(\MP_o^3;\tenset{3}{2})$, we have
\begin{align}
Z^-(s\hat T^2_2)_{12}=-Z^-(sT^2)_{13}=-s Z^{-}(T^2)_{13}
\end{align}
for any $s>0$.
Together with \eqref{eq604-4aaa}, we have
\begin{align}
\lambda^{2/3}Z^-(\lambda^{1/3}s\hat T^2_2)_{12}=-\lambda s Z^{-}(T^2)_{13}=  2 b(1;1) \lambda s = 2 b \lambda s.
\end{align}
In addition, both the permutations
$\ab{\begin{smallmatrix}1 & 2 & 3 \\ 1 & 3& 2 \end{smallmatrix}}$
and
$\ab{\begin{smallmatrix}1 & 2 & 3 \\ 2 & 3& 1 \end{smallmatrix}}$
send $\hat T^2_2$ to $T^2$.
Hence $Z^-(s\hat T^2_2)_{22}=-Z^-(s T^2)_{33}=Z^-(s T^2)_{33}$ which shows $Z^-(s\hat T^2_2)_{22}=0$.
Now \eqref{eq604-2a} and \eqref{eq604-4a} show
\begin{align}
&Z^-(sT^3)(e_1,e_2) + \lambda^{2/3}Z^-(\lambda^{1/3}s\hat T^2_2)\ab{\frac{1}{\lambda}e_1-\frac{1-\lambda}{\lambda}e_2,e_2}\\
&=cs^4+bs + \lambda^{2/3}Z^-(\lambda^{1/3}s\hat T^2_2)\ab{\frac{1}{\lambda}e_1,e_2}
=cs^4+3bs
\end{align}
and
\begin{align}
&\lambda^{2/3} Z^-(\lambda^{1/3}sT^3) \ab{\frac{1}{\lambda}e_1-\frac{1-\lambda}{\lambda}e_2,e_2} \\
&\qquad  +(1-\lambda)^{2/3}Z^-((1-\lambda)^{1/3}sT^3)\ab{e_1,-\frac{\lambda}{1-\lambda}e_1+\frac{1}{1-\lambda} e_2} \\
&= c \lambda s^4 + bs + c (1-\lambda) s^4 + bs = c s^4 + 2 bs.
\end{align}
Together with \eqref{30} for $Z^-$ with $n=3$ and $p=2$, we conclude $b=0$ and hence \eqref{eq:ind2a-} holds in this setting and \eqref{eq:indlow1} also holds.

Applying induction to $sT^1$ with  \eqref{eq:ind1a}, Theorems \ref{thm:p=1}, \ref{thm:p=n=2} and \ref{thm:contrav}, we get \eqref{eq:indlow2} for $p=2$.
\end{proof}

Next, we will use the above preliminary induction results to prove Theorem \ref{thm:low} in \S \ref{sec:low} and Theorem \ref{thm:tensor:po} and Theorem \ref{thm:tensor:neg} in \S \ref{sec:full}.

\subsection{Proof of Theorem \ref{thm:low}}\label{sec:low}
We prove Theorem \ref{thm:low} in this subsection and separate it into the following Lemmas \ref{Lem:Zo}, \ref{lem:low}.
\begin{lem}\label{Lem:Zo}
Let $n \ge p \ge 2$ and $Z \in \tenval(\MPon;\tenset{n}{p})$.
If $n>p$, then
\begin{align*}
Z(\{o\})= 0.
\end{align*}
If $n=p$, then there is a constant $c$ such that
\begin{align*}
Z(\{o\})= c \lct.
\end{align*}
\end{lem}
\begin{proof}
Define $\phi \in \sln$ by
\begin{align*}
\phi e_j:=r_je_j,~~ 1 \leq j \leq n,
\end{align*}
where $r_j>0$ and $r_{1} \cdots r_{n}=1$.
By the $\sln$ contravariance of $Z$,
\begin{align*}
Z(\{o\})_{\a_1 \dots \a_p}
=Z(\phi^{-t} \{o\})_{\a_1 \dots \a_p}
=(r_{1})^{\gamma_{1}} \cdots (r_{n})^{\gamma_{n}} Z(\{o\})_{\a_1 \dots \a_p},
\end{align*}
for every $\a_1, \dots, \a_p \in \{1,\dots,n\}$,
where $\gamma_{j}$ denotes how many times $j$ appears in $\{\a_1, \dots, \a_p\}$.

If $\seqalpha{\gamma}{1}{n}$ are not the same, then by the arbitrariness of $r_{1}, \dots, r_{n}$, we have
\begin{align*}
Z(\{o\})_{\a_1 \dots \a_p}=0
\end{align*}
for every $\a_1, \dots, \a_p \in \{1,\dots,n\}$.
But since $n \ge p$, if $\seqalpha{\gamma}{1}{n}$ are the same, then $n=p$ and $\seqalpha{\a}{1}{n}$ must be pairwise distinct.

Now set $c=Z(\{o\})_{1 \dots n}$. We only need to show that $Z(\{o\})_{\a_1 \dots \a_n} = c \lct_{\a_1 \dots \a_n}$ when $\seqalpha{\a}{1}{n}$ are pairwise distinct.
In fact, define $\psi \in \sln$ by $\psi e_1:=\lct_{\a_1 \dots \a_n} e_{\a_1}$ and $\psi e_i:=e_{\a_i}$ for all $i=\seqnb{2}{n}$.
Then we have
\begin{align*}
Z(\{o\})(e_{\a_1}, \dots , e_{\a_n})
&=Z(\psi \{o\})(e_{\a_1}, \dots , e_{\a_n})\\
&=Z(\{o\})(\psi^{-1}e_{\a_1}, \dots , \psi^{-1}e_{\a_n})
=\lct_{\a_1 \dots \a_n} Z(\{o\})(e_{1}, \dots , e_{n}).
\end{align*}
\end{proof}

\begin{lem}\label{lem:low}
Let $n \ge 3$ and $p \ge 2$ such that $n \ge p$, and let $Z \in \tenval(\MPon;\tenset{n}{p})$.
For $n>p$, $Z$ is simple.
For $n=p$, there is a constant $c_0$ such that
\begin{align*}
Z(sT^d)= c_0 \lct
\end{align*}
for all $1 \le d \le n-1$ and $s>0$.
\end{lem}

\begin{proof}
The proof is divided into three cases. Let $\seqalpha{\a}{1}{p} \in \seq{1}{n}$. The first case assumes that at least one of $\alpha_k$ for $1 \le k \le p$ lies in $\{1,\dots,d\}$. The second and third cases both assume that all $\alpha_k$ lie in $\{d+1,\dots,n\}$ with the second case for $d \ge 2$ and the third case for $d = 1$.

\textbf{Case 1.} There are $i \in \seq{1}{d}$, $q \in \seq{1}{p}$ and $1 \le k_1 < \dots < k_q \le p$ such that $\a_{k_1} =\a_{k_2}=\dots =\a_{k_q}=i$, while the other $\a_l$ are different from $i$. When $q=1$, we briefly write $k$ instead of $k_1$.

Let $s>0$. By \eqref{eq:indlow1} and \eqref{eq:indlow2}, 
there are $c_0'(\indfconst{i}{k_1,\dots,k_q}), c_0(\indfconst{i}{k_1,\dots,k_q}) \in \R$ which depend only on the indices $i$, $k_1,\dots,k_q$ such that
\begin{equation}\label{eq509-1}
\begin{aligned}
&Z(sT^d)_{\a_1 \dots \a_p}=\begin{cases}
c_0(\indfconst{i}{k}) \lct^{e_i^\bot}_{\a_1 \dots \a_{\hat{k}} \dots \a_p}, & \text{~if~} q=1 \text{~and~} p=n,\\
c_0(\indfconst{i}{1,\dots,p})s^{-p}, & \text{~if~} q=p,\\
0, & \text{~otherwise}.
\end{cases}
\end{aligned}
\end{equation}
for $d \ge 2$ and
\begin{equation}\label{eq509-2}
\begin{aligned}
&Z(sT^1)_{\a_1 \dots \a_p}=\begin{cases}
\ab{c_0'(\indfconst{1}{k})+c_0(\indfconst{1}{k})} \lct^{e_1^\bot}_{\a_1 \dots \a_{\hat{k}} \dots \a_p}, & \text{~if~} q=1 \text{~and~} p=n,\\
\ab{c_0'(\indfconst{1}{1,\dots,p})+c_0(\indfconst{1}{1,\dots,p})}s^{-p}, & \text{~if~} q=p,\\
0, & \text{~otherwise}.
\end{cases}
\end{aligned}
\end{equation}

If $p=n$ and $q=1$, to conclude that
\begin{align*}
Z(sT^d)_{\a_1 \dots \a_n}= c_0 \lct_{\a_1 \dots \a_n}
\end{align*}
for some constant $c_0 \in \R$ from \eqref{eq509-1} and \eqref{eq509-2}, we only need to assume that $\seqalpha{\a}{1}{n}$ are pairwise distinct numbers and show that \\
(i) $Z(T^d)_{\a_1 \dots \a_n}$ changes the sign when we switch any two indices $\a_k$ and $\a_l$ for all $d \ge 2$; \\
(ii) $Z(T^2)_{\a_1 \dots \a_n}=Z(T^1)_{\a_1 \dots \a_n}$.

For (i), since $q=1$, there is a $k$ such that $\a_k=1$, by \eqref{eq509-1}, we have
\begin{align*}
Z(T^d)_{\relseq{\a}{1}{n}{}} = c_0(\indf{1}{k}{1}) \lct_{\a_1 \dots \a_{\hat{k}} \dots \a_n} = Z(T^{n-1})_{\relseq{\a}{1}{n}{}}
\end{align*}
for all $d \ge 2$.
Assume first $\a_1,\a_2 \neq n$. We show $Z(T^{n-1})_{\a_1\a_2\a_3\dots\a_n}=-Z(T^{n-1})_{\a_2\a_1\a_3\dots\a_n}$.
In fact, define $\psi \in \sln$ with $\psi e_{\a_1}:=e_{\a_2}$, $\psi e_{\a_2}:=e_{\a_1}$, $\psi e_n:=-e_n$ and $\psi e_j:=e_j$ for the other $j \in \seq{1}{n}$.
Since only one $n$ appears among $\a_3,\dots,\a_n$, we have
\begin{align*}
Z(T^{n-1})_{\a_1\a_2\a_3\dots\a_n}=Z(\psi T^{n-1})_{\a_1\a_2\a_3\dots\a_n}=-Z( T^{n-1})_{\a_2\a_1\a_3\dots\a_n}.
\end{align*}
Now assume only one of $\a_1$ and $\a_2$ is $n$, for instance, $\a_1=1$, $\a_2=n$ and $\a_3=2$. We also want to show $Z(T^{n-1})_{\a_1\a_2\a_3\dots\a_n}=-Z(T^{n-1})_{\a_2\a_1\a_3\dots\a_n}$.
By \eqref{eq509-1}, we have
\begin{align*}
Z(T^{n-1})_{\a_1\a_2\a_3\a_4\dots\a_n}
=c_0(\indf{2}{3}{1}) \lct^{e_2^\bot}_{\a_1\a_2\a_4\dots\a_n}
&=-c_0(\indf{2}{3}{1}) \lct^{e_2^\bot}_{\a_2\a_1\a_4\dots\a_n}\\
&=-Z(T^{n-1})_{\a_2\a_1\a_3\a_4\dots\a_n}.
\end{align*}
Similarly, $Z(T^{n-1})_{\a_1 \dots \a_n}$ changes the sign when we switch any two indices $\a_k$ and $\a_l$.

For (ii), by \eqref{30a}, we have
\begin{align*}
&Z (T^2)(e_1, e_2, e_3, \dots, e_{n-1},e_n) + Z (T^1)\ab{\frac{1}{\lambda} e_1 - \frac{1-\lambda}{\lambda} e_2, e_2, e_3, \dots, e_{n-1}, \lambda e_n} \\
&=Z (T^2)\ab{\frac{1}{\lambda} e_1 - \frac{1-\lambda}{\lambda} e_2, e_2, e_3, \dots, e_{n-1}, \lambda e_n} \\
&\qquad + Z(T^2)\ab{e_1,-\frac{\lambda}{1-\lambda}e_1+\frac{1}{1-\lambda}e_2, e_3, \dots, e_{n-1}, (1-\lambda)e_n}
\end{align*}
for every $s >0$ and $\lambda \in (0,1)$.
Since those tensors are $n$-linear functionals, together with \eqref{eq509-1} and \eqref{eq509-2}, we have
\begin{align}
Z (T^2)(e_1, e_2, e_3, \dots, e_{n-1},e_n) = Z (T^1)(e_1, e_2, e_3, \dots, e_{n-1},e_n).
\end{align}
The statement (ii) for general pairwise distinct $\a_1,\dots,\a_n$ is similar.

Now we need to prove that $c_0(\indfconst{i}{1,\dots,p})= c_0'(\indfconst{1}{1,\dots,p})=0$ in \eqref{eq509-1} and \eqref{eq509-2} for $p=q$ to complete Case 1.
Applying the $\sln$ contravariance of $Z$, we easily find that $c_0(\indfconst{i}{1,\dots,p})=c_0(\indfconst{1}{1,\dots,p})$ for all $i \in \seq{1}{d}$.
Let $\lambda \in (0,1)$. The relation \eqref{30a} gives
\begin{align*}
&Z (T^2)(e_2, \dots, e_2) + Z (T^1)\ab{e_2, \dots, e_2} \\
=&Z (T^2)\ab{e_2, \dots, e_2}+ Z(T^2)\ab{-\frac{\lambda}{1-\lambda}e_1+\frac{1}{1-\lambda}e_2, \dots, -\frac{\lambda}{1-\lambda}e_1+\frac{1}{1-\lambda}e_2}.
\end{align*}
Together with \eqref{eq509-1}, we obtain
\begin{align*}
Z (T^1)\ab{e_2, \dots, e_2}
=\ab{\ab{-\frac{\lambda}{1-\lambda}}^p+ \ab{\frac{1}{1-\lambda}}^p} c_0(\indfconst{1}{1,\dots,p}).
\end{align*}
The left-hand side does not depend on $\lambda$, while the right-hand side does (since $p \ge 2$). Thus
$$c_0(\indfconst{1}{1,\dots,p})=Z (T^1)\ab{e_2, \dots, e_2}=0$$
and hence $c_0(\indfconst{1}{1,\dots,p})=0$ for all $i \in \seq{1}{d}$.
The relation \eqref{30a} also gives
\begin{align*}
&Z (T^2)(e_1, \dots, e_1) + Z (T^1)\ab{\frac{1}{\lambda} e_1 - \frac{1-\lambda}{\lambda} e_2, \dots, \frac{1}{\lambda} e_1 - \frac{1-\lambda}{\lambda} e_2} \\
=&Z (T^2)\ab{\frac{1}{\lambda} e_1 - \frac{1-\lambda}{\lambda} e_2, \dots, \frac{1}{\lambda} e_1 - \frac{1-\lambda}{\lambda} e_2}+ Z(T^2)\ab{e_2,\dots,e_2}.
\end{align*}
Together with \eqref{eq509-1} and \eqref{eq509-2}, we have
\begin{align}
\frac{1}{\lambda^p} c_0'(\indfconst{1}{1,\dots,p})= \ab{\ab{-\frac{\lambda}{1-\lambda}}^p+ \ab{\frac{1}{1-\lambda}}^p} c_0(\indfconst{1}{1,\dots,p}) =0.
\end{align}
That completes the proof of Case 1.

\textbf{Case 2.}
Let $d \ge 2$ and $\a_1, \dots, \a_p \in \{d+1,\dots,n\}$.
By the $\sln$ contravariance of $Z$, we may assume that $n \in \set{\seqalpha{\a}{1}{p}}$.
Let $\lambda \in (0,1)$. By \eqref{30a},
\begin{align*}
\lambda^{\gamma} Z(s\hat{T}^{d-1}_2)_{\a_1 \dots \a_p} = (\lambda^{\gamma} + (1-\lambda)^{\gamma} -1) Z(sT^d)_{\a_1 \dots \a_p},
\end{align*}
where $\gamma$ denotes how many times $n$ appears in $\{\a_1, \dots, \a_p\}$.
Note that $\gamma \ge 1$.
Letting $\lambda \to 1$ yields
\begin{align}\label{eq408-2}
Z(s\hat{T}^{d-1}_2)_{\a_1 \dots \a_p}=0
\end{align}
for all $\a_1 \dots \a_p \in \{d+1,\dots,n\}$.
Then
\begin{align}\label{eq607-1}
Z(sT^d)_{\a_1 \dots \a_p}=0
\end{align}
for all $\a_1, \dots, \a_p \in \seq{d+1}{n}$ with $d \ge 2$ and at least two of $\a_1, \dots, \a_p$ are the same ($\gamma>1$).

Next, let $\seqalpha{\a}{1}{p}$ be pairwise distinct in $\seq{d+1}{n}$ and $d \ge 2$. Observe that now $p$ must be smaller than or equal to $n-d$.
Since $p \ge 2$, we have $d \le n-2$.

Choose $\phi_3, \phi_4 \in \sln$ such that
\begin{align*}
\phi_3(e_1)=e_2,~ \phi_3(e_2)=e_1, \phi_3(e_{\a_1})=e_{\a_2},~ \phi_3(e_{\a_2})=e_{\a_1}, ~\phi_3(e_j)=e_j \text{~for~all~other~} j,
\end{align*}
and
\begin{align*}
\phi_4(e_{\a_1})=e_{\a_2},~ \phi_4(e_{\a_2})=-e_{\a_1}, ~\phi_4(e_j)=e_j \text{~for~all~other~} j.
\end{align*}
Then
\begin{align*}
Z(sT^d)_{\a_1\a_2\a_3\dots\a_p}=Z(\phi_3sT^d)_{\a_1\a_2\a_3\dots\a_p}=(\phi_3^{-t} \cdot Z(sT^d))_{\a_1\a_2\a_3\dots\a_p}=Z(sT^d)_{\a_2\a_1\a_3\dots\a_p}
\end{align*}
and
\begin{align*}
Z(sT^d)_{\a_1\a_2\a_3\dots\a_p}=Z(\phi_4sT^d)_{\a_1\a_2\a_3\dots\a_p}=(\phi_4^{-t} \cdot Z(sT^d))_{\a_1\a_2\a_3\dots\a_p}=-Z(sT^d)_{\a_2\a_1\a_3\dots\a_p}.
\end{align*}
Thus $Z(sT^d)_{\a_1\dots\a_p}=0$ for this case.

\textbf{Case 3.} Let $d =1$ and $\a_1 \dots \a_p \in \{2,\dots,n\}$.
We want to prove
\begin{equation}\label{eqn:st13}
Z(s T^{1})_{\a_1 \dots \a_p}=0.
\end{equation}

By \eqref{eq408-2} for $d=2$, \eqref{eqn:st13} holds for all $\a_1, \dots \a_p \in \seq{3}{n}$.

For $\seqalpha{\a}{1}{p} \in \seq{2}{n}$, assume w.l.o.g. $\a_1=\dots=\a_m=2$ and $\seqalpha{\a}{m+1}{p} \in \seq{3}{n}$ where $m\in\seq{1}{p}$.
Let $\lambda \in (0,1)$. Applying \eqref{30a} again, we have
\begin{align*}
&Z(sT^2)(e_{2},\dots,e_2,  e_{\a_{m+1}}, \dots, e_{\a_p}) + \lambda^{\gamma} Z(s T^{1})(e_{2},\dots,e_2,  e_{\a_{m+1}}, \dots, e_{\a_p}) \\
=&\lambda^{\gamma} Z(sT^2)(e_{2},\dots,e_2,  e_{\a_{m+1}}, \dots, e_{\a_p})
\\
&+ (1-\lambda)^{\gamma} Z(sT^2) \ab{-\frac{\lambda}{1-\lambda}e_1+\frac{1}{1-\lambda}e_2,\dots,-\frac{\lambda}{1-\lambda}e_1+\frac{1}{1-\lambda}e_2,  e_{\a_{m+1}}, \dots, e_{\a_p}}.
\end{align*}
As above, $\gamma$ denotes how many times $n$ appears in $\{\a_1, \dots, \a_p\}$.
Due to Case 1,
\begin{align*}
&Z(sT^2)(e_{2},\dots,e_2,  e_{\a_{m+1}}, \dots, e_{\a_p})\\
=&Z(sT^2) \ab{-\frac{\lambda}{1-\lambda}e_1+\frac{1}{1-\lambda}e_2,\dots,-\frac{\lambda}{1-\lambda}e_1+\frac{1}{1-\lambda}e_2,  e_{\a_{m+1}}, \dots, e_{\a_p}}=0
\end{align*}
if $n>p$.
Thus \eqref{eqn:st13} holds for this setting.
Now, we assume $n=p$. If $m=1$, there must be at least two of $\a_2,\dots,\a_p \in \seq{3}{n}$ coincide,
which, by Case 1, forces
\[Z(sT^2)_{2\a_2\dots\a_p}=Z(sT^2)_{1\a_2\dots\a_p}=0,\]
and then \eqref{eqn:st13} follows.
If $m=2$, by Case 1, we have
\[Z(sT^2)_{22\a_3\dots\a_p}=Z(sT^2)_{11\a_3\dots\a_p}=0,\]
and
\[Z(sT^2)_{12\a_3\dots\a_p}=-Z(sT^2)_{21\a_3\dots\a_p},\]
which together imply \eqref{eqn:st13}.
If $m>2$, Case 1 also shows
$$Z(sT^2) \ab{-\frac{\lambda}{1-\lambda}e_1+\frac{1}{1-\lambda}e_2,\dots,-\frac{\lambda}{1-\lambda}e_1+\frac{1}{1-\lambda}e_2,  e_{\a_{m+1}}, \dots, e_{\a_p}}=0.$$
Then we obtain \eqref{eqn:st13}, which completes this case.
\end{proof}

We can now prove Theorem \ref{thm:low}.
\begin{proof}[Proof of Theorem \ref{thm:low}.]
The proof follows directly from Lemmas \ref{Lem:Zo} and \ref{lem:low} with setting $c_0':=c-c_0$ for constants $c,c_0$ in Lemmas \ref{Lem:Zo} and \ref{lem:low}, respectively.
\end{proof}

\subsection{Proofs of Theorem \ref{thm:tensor:po} and Theorem \ref{thm:tensor:neg}}\label{sec:full}
We always assume that $Z^+ \in \tenval_0(\MPon;\tenset{n}{p})$ and $Z^- \in \tenval_1(\MPon;\tenset{n}{p})$ are simple valuations in this subsection. Then the induction step also works for $Z^+$ and $Z^-$ and \eqref{eq:ind2}, \eqref{eq:ind2a}, \eqref{eq:ind2a-}, \eqref{eq:ind2b-}, \eqref{eq:ind2bbb-}, \eqref{eq:ind2c-} and \eqref{eq:ind2d-} also hold for the simple $Z^+$ and $Z^-$ with $c_0 \equiv 0$ appeared in \eqref{eq:ind2a}, \eqref{eq:ind2b-}, \eqref{eq:ind2bbb-} and \eqref{eq:ind2d-}.
Based on these, we further study the coefficients in \eqref{eq:ind2}, \eqref{eq:ind2a}, \eqref{eq:ind2a-}, \eqref{eq:ind2b-}, \eqref{eq:ind2bbb-}, and \eqref{eq:ind2d-} in Lemmas \ref{lem:coef2}, \ref{lem:coef3} and \ref{lem:coef4}. Then, after establishing Lemma \ref{lem:coef5}, we complete the proofs of Theorem \ref{thm:tensor:po} and Theorem \ref{thm:tensor:neg}.

\begin{lem}\label{lem:coef2}
The function $\zeta$ in \eqref{eq:ind2} and \eqref{eq:ind2a} has the following properties:
\begin{equation}
	\zeta(\indf{i}{k_1,\dots,k_{q}}{s^n})=\zeta(\indf{1}{1}{s^n})
\end{equation}
for any $1 \leq q \leq p$, $i \in \{1,\dots,n\}$, $1 \leq k_1 < \dots < k_{q} \leq p$ and $s>0$.
\end{lem}

\begin{proof}
First, let $q=p$ and assume $i \neq 1$.
We can choose a $\phi \in \slpm{n}$ such that $\phi$ switches between $e_i$ and $e_1$ while fixing the other $e_j$ for $j \notin \{1,i\}$.
By the $\slpm{n}$-$0$-contravariance of $Z^+$,
\eqref{eq:ind2} and \eqref{eq:ind2a} imply that
\begin{align*}
\zeta(\indf{i}{1,\dots,p}{s^n})s^{-p}=Z^+(sT^n)_{i\dots i}\
&=Z^+(sT^n)_{j\dots j}\\
&=\zeta(\indf{j}{1,\dots,p}{s^n})s^{-p}
\end{align*}
for any $i,j \in \{1,\dots,n\}$, which confirms
\begin{align}\label{eq:zeta1}
\zeta(\indf{i}{1,\dots,p}{s^n})
=\zeta(\indf{1}{1,\dots,p}{s^{n}})
\end{align}
for any $i \in \{1,\dots,n\}$.

Next, let $q \le p-1$ and assume $i \neq j$.
Set $\a_1 \dots \a_{\hat{k}_1} \dots \a_{\hat{k}_q} \dots \a_p$ to be all pairwise distinct numbers with $\a_l=j$. By \eqref{eq:ind2} and \eqref{eq:ind2a}, we get
\begin{equation}\label{eq617-1}
\begin{aligned}
\zeta(\indf{i}{k_1,\dots,k_q}{s^n})s^{-p} &=Z^+(sT^n)_{\a_1 \dots \a_{k_1} \dots \a_{k_q} \dots \a_p}\\
&=Z^+(sT^n)_{\a_1 \dots \a_{l} \dots \a_p}
=\zeta(\indf{j}{l}{s^n}) s^{-p} .
\end{aligned}
\end{equation}
Now we prove that
\begin{align}\label{eq616-2}
\zeta(\indf{i}{k}{s^n}) = \zeta(\indf{j}{l}{s^n})
\end{align}
for every $i,j \in \seq{1}{n}$ and $k,l\in\seq{1}{p}$.
Assume first $i \neq j$ and $k=l$. Let $\a_1 \dots \a_{\hat{k}} \dots \a_p \in \seq{1}{n} \setminus \set{i,j}$.
We still apply the $\slpm{n}$-$0$-contravariance of $Z^+$ similar as in the case $q=p$ to get
\begin{equation}\label{eq616-1}
\begin{aligned}
\zeta(\indf{i}{k}{s^n}) s^{-p}&= Z^+(sT^n)_{\a_1 \dots \a_{k-1} \a_k \a_{k+1}  \dots \a_p} \\
&= Z^+(sT^n)_{\a_1 \dots \a_{k-1} \beta_{k} \a_{k+1} \dots \a_p} = \zeta(\indf{j}{k}{s^n})s^{-p},
\end{aligned}
\end{equation}
where $\a_k=i$ and $\beta_k=j$.
Second, assume $k \neq l$ and $i \neq j$. Set $\a_k=i$, $\a_l=j$ and $\a_1, \dots, \a_{\hat{k}}, \dots, \a_{\hat{l}} ,\dots, \a_p \in \{1,\dots, n\} \setminus \{i,j\}$ in \eqref{eq:ind2} and \eqref{eq:ind2a}. We have
\begin{align*}
\zeta(\indf{i}{k}{s^n}) s^{-p} &= Z^+(sT^n)_{\a_1 \dots \a_{k} \dots \a_l \dots \a_p}
= \zeta(\indf{j}{l}{s^n}) s^{-p}.
\end{align*}
Together with \eqref{eq616-1}, it turns out
\begin{align*}
\zeta(\indf{i}{k}{s^n}) = \zeta(\indf{j}{l}{s^n})
= \zeta(\indf{i}{l}{s^n}).
\end{align*}
Considering the above equalities covers all cases (indices both different, first different/second same, first same/second different), so \eqref{eq616-2} follows for all $i,j \in \seq{1}{n}$ and $k,l\in\seq{1}{p}$.
With \eqref{eq617-1}, we obtain
\begin{align}\label{eq:zeta2}
\zeta(\indf{i}{k_1,\dots,k_{q}}{s^n})
=\zeta(\indf{1}{1}{s^n})
\end{align}
for any $1 \leq q \leq p-1$, $i \in \{1,\dots,n\}$, and $1 \leq k_1 \le \dots \le k_{q} \leq p$.

Finally, by \eqref{30} for the simple valuation $Z^+$ together with \eqref{eq:ind2}, \eqref{eq:ind2a}, \eqref{eq:zeta1} and \eqref{eq:zeta2}, we obtain
\begin{align*}
&\zeta(\indf{1}{1,\dots,p}{s})s^{-p/n}=Z^+(s^{1/n}T^n)_{1\dots 1} \\
=&\lambda^{p/n} Z^+((\lambda s)^{1/n}T^n) \ab{\frac{1}{\lambda}e_1-\frac{1-\lambda}{\lambda}e_2, \cdots, \frac{1}{\lambda}e_1-\frac{1-\lambda}{\lambda}e_2,}\\
&+(1-\lambda)^{p/n}Z^+ \big(((1-\lambda)s)^{1/n}T^n \big)\ab{e_1,\dots,e_1}\\
=&\lambda^{-p}s^{-p/n} \ab{\zeta(\indf{1}{1,\dots,p}{\lambda s})(1+(-1)^p(1-\lambda)^p)+ \sum_{j=1}^{p-1} \binom{p}{j} (-1)^j (1-\lambda)^j\zeta(\indf{1}{1}{\lambda s})} \\
&+ s^{-p/n}\zeta(\indf{1}{1,\dots,p}{(1-\lambda) s}).
\end{align*}
Then the additive property of Cauchy functions yields
\begin{align*}
&\ab{1-\lambda^{-p}\ab{1+(-1)^p(1-\lambda)^p}}\zeta(\indf{1}{1,\dots,p}{\lambda s}) \\
=&\sum_{j=1}^{p-1} \lambda^{-p}\binom{p}{j} (-1)^j (1-\lambda)^j \zeta(\indf{1}{1}{\lambda s})
\end{align*}
for any $\lambda s>0$.
Clearly
\begin{align*}
1-\lambda^{-p}\ab{1+(-1)^p(1-\lambda)^p}=\sum_{j=1}^{p-1} \lambda^{-p} \binom{p}{j} (-1)^j (1-\lambda)^j.
\end{align*}
Hence
\begin{equation}
\begin{aligned}
\zeta(\indf{1}{1,\dots,p}{s^{n}})
=\zeta(\indf{1}{1}{s^n}),
\end{aligned}
\end{equation}
and the desired result follows from \eqref{eq:zeta1} and \eqref{eq:zeta2}.
\end{proof}

\begin{lem}\label{lem:coef3}
The constant $c(\indfconst{i}{k})$ in \eqref{eq:ind2a-} for the case $p=n-1$ satisfies the following properties:
\begin{align}\label{eq:tau1}
(-1)^{i+k}c(\indfconst{i}{k})
=c(\indfconst{1}{1})
\end{align}
for any $i \in \{1,\dots,n\}$ and $k \in \{1,\dots,p\}$.
\end{lem}
\begin{proof}
Assume first $i \neq j$ and $k=l$. Choose $\alpha_1 \dots \alpha_{\hat{k}} \dots \alpha_{n-1}$ to be pairwise distinct in $\seq{1}{n} \setminus \{i,j\}$.
By \eqref{eq:ind2a-} for $p=n-1$ and $Z^- \in \tenval_1(\MPon;\tenset{n}{p})$, we have
\begin{align*}
c(\indfconst{i}{k})s^{n+1} \lct_{\a_1 \dots \a_{\hat{k}} \dots \a_{n-1}}^{e_i^\bot} &=Z^-(sT^n)_{\a_1 \dots \a_{k}\dots \a_{n-1}} \\
&= -Z^-(sT^n)_{\a_1 \dots \beta_{k}\dots \a_{n-1}}
=- c(\indfconst{j}{k})s^{n+1} \lct_{\a_1 \dots \a_{\hat{k}} \dots \a_{n-1}}^{e_j^\bot}
\end{align*}
with $\a_k=i$ and $\beta_k=j$.
Note that
\begin{align*}
\lct_{\a_1 \dots \a_{\hat{k}} \dots \a_{n-1}}^{e_i^\bot}
&= \det \ab{e_{\a_1},\dots,e_{\a_{k-1}},e_{\a_{k+1}},\dots,e_{\a_{n-1}},e_j}_{(n-1)} \\
&= (-1)^{i+k} \det \ab{e_{\a_1},\dots,e_{\a_{k-1}},e_i,e_{\a_{k+1}},\dots,e_{\a_{n-1}},e_j} \\
&= (-1)^{i+k+1} \det \ab{e_{\a_1},\dots,e_{\a_{k-1}},e_j,e_{\a_{k+1}},\dots,e_{\a_{n-1}},e_i} \\
&=(-1)^{i-j+1} \det \ab{e_{\a_1},\dots,e_{\a_{k-1}},e_{\a_{k+1}},\dots,e_{\a_{n-1}},e_i}_{(n-1)}\\
&=(-1)^{i-j+1}\lct_{\a_1 \dots \a_{\hat{k}} \dots \a_{n-1}}^{e_j^\bot}.
\end{align*}
Hence
\begin{align}\label{eq618-1}
c(\indfconst{i}{k})=(-1)^{i-j}c(\indfconst{j}{k}).
\end{align}

Second, assume $i\neq j$ and $k \neq l$. Set $\a_k=i$, $\a_l=j$ and let $\a_1, \dots, \a_{\hat{k}}, \dots, \a_{\hat{l}} ,\dots \a_{n-1}$ be pairwise distinct in $\{1,\dots, n\} \setminus \{i,j\}$ in \eqref{eq:ind2a-}. We have
\begin{align*}
c(\indfconst{i}{k})s^{n+1} \lct_{\a_1 \dots \a_{\hat{k}} \dots \a_{n-1}}^{e_i^\bot} &=Z^-(sT^n)_{\a_1 \dots \a_{k} \dots \a_{l} \dots \a_{n-1}} \\
&= c(\indfconst{j}{l})s^{n+1} \lct_{\a_1 \dots \a_{\hat{l}} \dots \a_{n-1}}^{e_j^\bot},
\end{align*}
and note that
\begin{align*}
\lct_{\a_1 \dots \a_{\hat{k}} \dots \a_{n-1}}^{e_i^\bot}
&= (-1)^{i+k} \det \ab{e_{\a_1},\dots,e_{\a_{k-1}},e_{i},e_{\a_{k+1}},\dots,e_{\a_{n-1}},e_{\a_n}} \\
&= (-1)^{i+k} \det \ab{e_{\a_1},\dots,e_{\a_{k}},\dots,e_{\a_{l}},\dots,e_{\a_{n-1}},e_{\a_n}} \\
&= (-1)^{i+k-j-l}\lct_{\a_1 \dots \a_{\hat{l}} \dots \a_{n-1}}^{e_j^\bot},
\end{align*}
with $\a_n \in \seq{1}{n} \setminus \seq{\a_1}{\a_{n-1}}$.
Hence
\begin{align*}
c(\indfconst{i}{k})=(-1)^{i+k-j-l} c(\indfconst{j}{l}).
\end{align*}
Now together with \eqref{eq618-1}, we obtain
\begin{align*}
c(\indfconst{i}{k})=(-1)^{i+k-j-l} c(\indfconst{j}{l})
=(-1)^{k-l} c(\indfconst{i}{l}).
\end{align*}
All together, we get
\begin{align*}
c(\indfconst{i}{k})=(-1)^{i+k-j-l} c(\indfconst{j}{l})=(-1)^{i+k} c(\indfconst{1}{1})
\end{align*}
for all $i \in \seq{1}{n}$ and $k \in \seq{1}{n-1}$.
\end{proof}

\begin{lem}\label{lem:coef4}
The constant $c_{r+1}(\indfconst{i}{k})$ and the function $\eta(\indf{i}{k}{s^n})$ in \eqref{eq:ind2b-}, \eqref{eq:ind2bbb-} and \eqref{eq:ind2d-} have the following property:
\begin{align*}
&(-1)^{i+k} \ab{\sum_{r=0}^{n-3} (-1)^{nr} c_{r+1}(\indfconst{i}{k})s^n +\eta(\indf{i}{k}{s^n})}\\
=&\sum_{r=0}^{n-3} (-1)^{nr} c_{r+1}(\indfconst{1}{1})s^{n} +\eta(\indf{1}{1}{s^{n}})
\end{align*}
for all $i,k \in \seq{1}{n}$ and $s>0$.
\end{lem}
\begin{proof}
Let $\a_1, \dots, \a_n$ be pairwise distinct in $\{1,\dots, n\}$ with $\a_{k}=i$.
Then \eqref{eq:ind2b-} and \eqref{eq:ind2bbb-} are simplified to \eqref{eq:ind2d-}.
Assume $i \neq j$ and $k \neq l$.
Apply the $\slpm{n}$ transform that switches $e_i$ and $e_j$ and fixes the other $e_m$ for $m \neq i,j$, to $sT^n$.
The $\slpm{n}$-$1$-contravariance of $Z^-$ together with \eqref{eq:ind2d-} gives
\begin{align*}
&(-1)^{i+k} \ab{\sum_{r=0}^{n-3} (-1)^{nr} c_{r+1}(\indfconst{i}{k})s^n +\eta(\indf{i}{k}{s^n})} \lct_{\a_1 \dots \a_{k-1} i \a_{k+1} \dots \a_{l-1} j \a_{l+1} \dots \a_n}^{(n)} \\
=&Z^-(sT^n)_{\a_1 \dots \a_{k-1} i \a_{k+1} \dots \a_{l-1} j \a_{l+1} \dots \a_n} \\
=&-Z^-(sT^n)_{\a_1 \dots \a_{k-1} j \a_{k+1} \dots \a_{l-1} i \a_{l+1} \dots \a_n}\\
=&-(-1)^{i+l} \ab{\sum_{r=0}^{n-3} (-1)^{nr} c_{r+1}(\indfconst{i}{l})s^n +\eta(\indf{i}{l}{s^n})} \lct_{\a_1 \dots \a_{k-1} j \a_{k+1} \dots \a_{l-1} i \a_{l+1} \dots \a_n}^{(n)} \\
=&-(-1)^{j+k} \ab{\sum_{r=0}^{n-3} (-1)^{nr} c_{r+1}(\indfconst{j}{k})s^n +\eta(\indf{j}{k}{s^n})} \lct_{\a_1 \dots \a_{k-1} j \a_{k+1} \dots \a_{l-1} i \a_{l+1} \dots \a_n}^{(n)}.
\end{align*}
Note that
\begin{align*}
\lct_{\a_1 \dots \a_{k-1} i \a_{k+1} \dots \a_{l-1} j \a_{l+1} \dots \a_n}^{(n)}
=-\lct_{\a_1 \dots \a_{k-1} j \a_{k+1} \dots \a_{l-1} i \a_{l+1} \dots \a_n}^{(n)}.
\end{align*}
Thus
\begin{align*}
&(-1)^{i+k} \ab{\sum_{r=0}^{n-3} (-1)^{nr} c_{r+1}(\indfconst{i}{k})s^n +\eta(\indf{i}{k}{s^n})} \\
=&(-1)^{i+l} \ab{\sum_{r=0}^{n-3} (-1)^{nr} c_{r+1}(\indfconst{i}{l})s^n +\eta(\indf{i}{l}{s^n})} \\
=&(-1)^{j+k} \ab{\sum_{r=0}^{n-3} (-1)^{nr} c_{r+1}(\indfconst{j}{k})s^n +\eta(\indf{j}{k}{s^n})}\\
=&(-1)^{j+l} \ab{\sum_{r=0}^{n-3} (-1)^{nr} c_{r+1}(\indfconst{j}{l})s^n +\eta(\indf{j}{l}{s^n})}.
\end{align*}
\end{proof}

\begin{lem}\label{lem:coef5}
Let $s>0$ and let $\a_1\dots\a_{k_1}\dots\a_{\hat k_2}\dots\a_l\dots\a_n$ be pairwise distinct in $\seq{1}{n}$ with $\a_{k_1}=\a_{k_2}$. Then
\begin{align*}
Z^-(sT^n)_{\a_1\dots\a_{k_1}\dots\a_{k_2}\dots\a_l\dots\a_n}
+Z^-(sT^n)_{\a_1\dots\a_{k_1}\dots\a_l\dots\a_{k_2}\dots\a_n}
+Z^-(sT^n)_{\a_1\dots\a_l\dots\a_{k_1}\dots\a_{k_2}\dots\a_n}=0.
\end{align*}
\end{lem}
\begin{proof}
We only prove the case in which $k_1=1$, $k_2=2$, $l=3$ with $\a_{1}=\a_{2}=1$ and $\a_3=2$. The other cases are analogous.

By \eqref{30} for the simple valuation $Z^-$, together with \eqref{eq:ind2c-}, we have
\begin{align*}
&Z^-(sT^n)_{112\a_4\dots\a_n}\\
=&\lambda Z^-(\lambda^{1/n}sT^n) \ab{\frac{1}{\lambda}e_1-\frac{1-\lambda}{\lambda}e_2,\frac{1}{\lambda}e_1-\frac{1-\lambda}{\lambda}e_2,e_2,e_{\a_4},\dots,e_{\a_n}}\\
&+(1-\lambda)Z^-((1-\lambda)^{1/n}sT^n)\ab{e_1,e_1,-\frac{\lambda}{1-\lambda}e_1+\frac{1}{1-\lambda}e_2,e_{\a_4},\dots,e_{\a_n}}\\
=&Z^-(sT^n)_{112\a_4\dots\a_n} -(1-\lambda)Z^-(sT^n)_{122\a_4\dots\a_n} -(1-\lambda)Z^-(sT^n)_{212\a_4\dots\a_n} \\
&+(1-\lambda)Z^-(sT^n)_{112\a_4\dots\a_n}.
\end{align*}
Thus
\begin{align*}
Z^-(sT^n)_{112\a_4\dots\a_n}-Z^-(sT^n)_{122\a_4\dots\a_n} -Z^-(sT^n)_{212\a_4\dots\a_n}=0.
\end{align*}
Further with the $\slpm{n}$-$1$-contravariance of $Z^-$, we conclude
\begin{align*}
Z^-(sT^n)_{112\a_4\dots\a_n}+Z^-(sT^n)_{211\a_4\dots\a_n}+Z^-(sT^n)_{121\a_4\dots\a_n}=0.
\end{align*}
\end{proof}

\begin{proof}[Proof of Theorem \ref{thm:tensor:po}.]
Recall that Lemma \ref{lem:coef2} shows that the function $\zeta(\indf{i}{k_1,\dots,k_{q}}{\cdot})$ in \eqref{eq:ind2} and \eqref{eq:ind2a} satisfies
$$\zeta(\indf{i}{k_1,\dots,k_{q}}{s^n}) s^{-p}=\zeta(\indf{1}{1}{s^n})s^{-p}$$
for all $s>0$, $1 \leq q \leq p$, $i \in \{1,\dots,n\}$ and $1 \leq k_1 <\dots <k_{q} \leq p$.

Set $\zeta(s^n)=\zeta(\indf{1}{1}{n! \, s^n})$ for $s>0$.
Clearly $\zeta$ is a Cauchy function on $(0,\infty)$.
Then \eqref{eq615-0} and \eqref{eq:ind2} and \eqref{eq:ind2a} imply the desired result.
\end{proof}

\begin{proof}[Proof of Theorem \ref{thm:tensor:neg}.]
The case $p \le n-2$ follows directly from \eqref{eq:ind2c-}.

For the case $p=n-1$, by \eqref{eq:ind2a-} and Lemma \ref{lem:coef3}, we have
\begin{align*}
Z^-(sT^n)_{\a_1 \dots \a_{k} \dots \a_{n-1}}
&=c(\indfconst{i}{k})s^{n+1} \lct_{\a_1 \dots \a_{\hat k} \dots \a_{n-1}}^{e_i^\bot}\\
&=c(\indfconst{i}{k})(-1)^{i+k} s^{n+1}\lct_{\a_1\dots \a_{k} \dots \a_{n-1}}^{(n)}
=c(\indfconst{1}{1})s^{n+1}\lct_{\a_1\dots \a_{k} \dots \a_{n-1}}^{(n)}
\end{align*}
for any pairwise distinct $\a_1,\dots,\a_{\hat k},\dots,\a_{n-1} \in \set{1,\dots,\hat i,\dots,n}$ with $\a_k=i$.
Set $c=(n+1)! \,c(\indfconst{1}{1})$. Then \eqref{eq615-1}, \eqref{eq615-2}, \eqref{eq:ind2a-} and \eqref{eq:ind2c-} imply the desired result.

For the case $p=n$, set
\begin{align*}
\tilde c_{n-k+1}:=(-1)^{n(n-k)} n! \, Z^{-}(T^n)_{\a_1\dots\a_{k}\dots\a_n}, \quad k \in \seq{2}{n},
\end{align*}
with $\a_1\dots\a_{\hat {k}}\dots\a_n=12\dots(n-1)$ and $\a_{k}=1$ and set
\begin{align*}
\tilde \eta\ab{\frac{1}{n!} s^n}:=\sum_{r=0}^{n-3} (-1)^{nr} c_{r+1}(\indfconst{1}{1})s^{n} +\eta(\indf{1}{1}{s^{n}})-\sum_{r=0}^{n-2} (-1)^{(n-1)r} \tilde c_{r+1} s^n
\end{align*}
for any $s>0$ and $\tilde \eta(0)=0$.
Clearly $\tilde \eta$ is a Cauchy function on $[0,\infty)$.
Define $\tilde Z:\MPon \to \tenset{n}{p}$ by
$$\tilde Z P := \sum_{r=0}^{n-2} \tilde c_{r+1} \ab{\sigma^{r} \left(\sum_{u \in \MN_o(sT^n)} \ct{m}{n}{F(sT^n,u)} \otimes u \right)} + \tilde \eta\big(V_n(sT^n)\big) \lct$$
for all $P \in \MPon$.
Theorem \ref{thm:contrav} shows that $\tilde Z \in \tenval_1(\MPon;\tenset{n}{p})$ and $\tilde Z$ is simple.
Now we only need to show that $\tilde Z(sT^n)=Z^-(sT^n)$ for any $s>0$; then, after an appropriate substitution of symbols, Theorem \ref{thm:tensor:neg} follows.

Direct calculation with the definition of the Levi-Civita symbol, \eqref{eq615-2}, \eqref{eq:ind2c-}, \eqref{eq:ind2d-} and Lemma \ref{lem:coef4} show that
\begin{align}\label{eq620-2}
Z^{-}(sT^n)_{\a_1\dots\a_n}&=\tilde Z(sT^n)_{\a_1\dots\a_n}
\end{align}
for all $\a_1,\dots,\a_n \in \seq{1}{n}$ except the remaining case where only two of $\a_1,\dots,\a_n$ are the same.
We now set $\a_1\dots\a_{\hat {k}}\dots\a_n=12\dots(n-1)$ and $\a_{k}=1$ for $k \in \seq{2}{n}$.
By \eqref{eq615-2} and $Z^{-}(sT^n)_{\a_1\dots\a_n}=s^nZ^{-}(T^n)_{\a_1\dots\a_n}$ (which follows from \eqref{eq:ind2b-} and \eqref{eq:ind2bbb-}),  we have
\begin{equation}\label{eq620-1}
\begin{aligned}
\tilde Z(sT^n)_{\a_1\dots\a_n}
&=\frac{1}{n!} s^n \sum_{r=0}^{n-2} c_{r+1} \lct_{\a_{\sigma^{-(r)}(1)}\dots\a_{\sigma^{-(r)}(n-1)}} 
+\eta\big(V_n(sT^n)\big) \lct_{\a_1\dots\a_n} \\
&=\frac{1}{n!} c_{n-k+1} s^n \lct_{k\dots (n-1) 1 \dots (k-1)n} \\
&=\frac{1}{n!} (-1)^{n(n-k)} c_{n-k+1} s^n
=Z^{-}(sT^n)_{\a_1\dots\a_{k}\dots\a_n}.
\end{aligned}
\end{equation}
That proves the desired result for this special setting.
For the general case in which $\beta_1,\dots,\beta_n \in \seq{1}{n}$ are with exactly one repeated element, Lemma \ref{lem:coef5} and the $\slpm{n}$-$1$-contravariance of $Z^-$ show that $Z^-(sT^n)_{\beta_1\dots\beta_n}$ is uniquely determined by $Z^-(sT^n)_{\a_1\dots\a_{k}\dots\a_n}$.
For instance, for $\beta_1=\beta_k=i$, define a permutation $\theta$ of $\seq{1}{n}$ with $\theta(1):=i$, $\theta(\a_l):=\beta_l$ for $l \in \seq{2}{n} \setminus \set{k}$ and $\theta(n):=j\in \seq{1}{n} \setminus \set{\beta_1,\dots,\beta_n}$.
Then the $\slpm{n}$-$1$-contravariance of $Z^-$ gives
$Z^-(sT^n)_{\beta_1\dots\beta_n}=(\sgn \theta)Z^-(sT^n)_{\a_1\dots\a_n}.$
Also, for $\beta_{k_1}=\beta_{k_2}=i$ with $1<k_1<k_2 \le n$, Lemma \ref{lem:coef5} with $l=1$ shows that this case is uniquely determined by the case that $\beta_1=\beta_{k_1}=i$ and $\beta_1=\beta_{k_2}=i$ which reduces to the previous situation.
Note that $\tilde Z \in \tenval_1(\MPon;\tenset{n}{p})$.
Then $\tilde Z(sT^n)_{\beta_1\dots\beta_n}$ is also uniquely determined by $\tilde Z(sT^n)_{\a_1\dots\a_n}$ with the same formulas.
Hence $\tilde Z(sT^n)_{\beta_1\dots\beta_n} = Z^{-} (sT^n)_{\beta_1\dots\beta_n}$ which completes the proof.
\end{proof}

\section{Extensions to \texorpdfstring{$\polyn$}{$\mathcal P^n$}}\label{sec:polyn}
In this section, we extend Theorems \ref{thm:tensor:p+2}, \ref{thm:tensor:p+1}, \ref{thm:tensor:p} to $\polyn$.
Recall that $[o,P]$ is the convex hull of $o$ and $P \in \MP^n$.

\begin{thm}\label{thm:1p}
	Let $n-2 \geq p \ge 2$.
	A mapping $Z: \polyn \to \tenset{n}{p}$ is an $\SLn$ contravariant valuation if and only if there are Cauchy functions $\zeta_1,\zeta_2: \R \to \R$ such that
	\begin{equation}\label{eqn:61}
		ZP= \sum_{u \in \mathcal {N}_o (P)} h_P^{-p}(u) \zeta_1(V(P,u)) \symtp{u}{p}
		+\sum_{u \in \mathcal {N}_o ([o,P])} h_{[o,P]}^{-p}(u) \zeta_2(V([o,P],u)) \symtp{u}{p}
	\end{equation}
	for every $P \in \polyn$.
\end{thm}
\begin{thm}\label{thm:2p}
	Let $n \ge 3$.
	A mapping $Z: \polyn \to \tenset{n}{n-1}$ is an $\SLn$ contravariant valuation if and only if there are Cauchy functions $\zeta_1,\zeta_2: \R \to \R$ and constants $c,\tilde c \in \R$ such that
	\begin{align*}
		ZP = &\sum_{u \in \mathcal {N}_o (P)}  h_P^{-(n-1)}(u) \zeta_1(V(P,u)) \symtp{u}{n-1}
		+ \sum_{u \in \mathcal {N}_o ([o,P])}  h_{[o,P]}^{-(n-1)}(u) \zeta_2(V([o,P],u)) \symtp{u}{n-1}\\
		&+ c \ct{m}{n+1}{P} + \tilde c \ct{m}{n+1}{[o,P]}
	\end{align*}
	for every $P \in \polyn$.
\end{thm}

\begin{thm}\label{thm:3p}
Let $n \ge 3$.
A mapping $Z: \polyn \to \tenset{n}{n}$ is an $\SLn$ contravariant valuation if and only if there are Cauchy functions $\zeta_1,\zeta_2,\eta_1,\eta_2: \R \to \R$ and constants $c_0',c_0, c_1,\dots,c_{n-1}, \tilde c_0,\tilde c_1,\dots,\tilde c_{n-1}$ such that
\begin{align*}
&ZP \\
&=\sum_{u \in \mathcal {N}_o (P)} h_P^{-n}(u) \zeta_1(V(P,u)) \symtp{u}{n}
+ \sum_{u \in \mathcal {N}_o ([o,P])} h_{[o,P]}^{-n}(u) \zeta_2(V([o,P],u)) \symtp{u}{n}\\
&~+ \sum_{r=0}^{n-2} c_{r+1} \sigma^r \left(\sum_{u \in \MN_o(P)} \ct{m}{n}{F(P,u)} \otimes u \right)
+ \sum_{r=0}^{n-2} \tilde c_{r+1} \sigma^r \left(\sum_{u \in \MN_o([o,P])} \ct{m}{n}{F([o,P],u)} \otimes u \right)\\
&~+ (c_0'(-1)^{\dim P}V_0(\{o\} \cap \relint P) + c_0V_0(P) +\tilde c_0V_0(\{o\} \cap P) + \eta_1\big(V_n(P)\big) + \eta_2(V_n([o,P]))) \lct
\end{align*}
for every $P \in \polyn$.

Let $n=2$.
A mapping $Z: \polye \to \tenset{2}{2}$ is an $\SL{2}$ contravariant valuation if and only if there are Cauchy functions $\zeta_1,\zeta_2,\eta_1,\eta_2: \R \to \R$ and constants $a,b,\tilde a,\tilde b$, $c_0',c_0, c_1, \tilde c_0,\tilde c_1$ such that
\begin{align*}
&ZP \\
& =\sum_{u \in \mathcal {N}_o (P)} h_P^{-2}(u) \zeta_1(V(P,u)) \symtp{u}{2}
+\sum_{u \in \mathcal {N}_o ([o,P])} h_{[o,P]}^{-2}(u) \zeta_2(V([o,P],u)) \symtp{u}{2}\\
&~+ c_{1}  \sum_{u \in \MN_o(P)} \ct{m}{2}{F(P,u) \otimes u }
+ \tilde c_1  \sum_{u \in \MN_o([o,P])} \ct{m}{2}{F([o,P],u) \otimes u}\\
&~+ (c_0'(-1)^{\dim P}V_0(\{o\} \cap \relint P) + c_0V_0(P) +\tilde c_0V_0(\{o\} \cap P) + \eta_1\big(V_n(P)\big) + \eta_2(V_n([o,P]))) \lct \\
&~+ a  M^{2,0}(\rho P) + \tilde a M^{2,0}(\rho[o,P]) + b \dten{2} (\rho P) +\tilde b  \dten{2} (\rho [o,P])
\end{align*}
for every $P \in \polye$.
\end{thm}

The proofs of Theorems \ref{thm:1p}, \ref{thm:2p} and \ref{thm:3p} are based on Lemma \ref{lemuq2} and are similar to the proofs of \cite[Theorem 2]{TLL2024}, \cite[Theorem 2]{LR2017sl} and \cite[Theorem 2]{SW12}.
We give the proof of Theorem \ref{thm:1p}, and omit others because
they follow by analogous arguments.
\begin{proof}[Proof of Theorem \ref{thm:1p}]
Let
\begin{align*}
	Z_1P=\sum_{u \in \MN_o (P)} h_P^{-p}(u) \zeta_1(V(P,u)) \symtp{u}{p} \text{ and }Z_2P=\sum_{u \in \MN_o ([o,P])} h_{[o,P]}^{-p}(u) \zeta_2(V([o,P],u)) \symtp{u}{p}
\end{align*}
for every \(P\in\MP^n\).
First, due to Theorem \ref{thm:contrav}, $Z_1$ is an $\sln$ contravariant valuation on $\MP^n$.
Next, for $P,Q\in\MP^n$ with $P\cup Q\in\MP^n$, we have $[o,P\cup Q]=[o,P]\cup[o,Q]$ and $[o,P\cap Q]=[o,P]\cap[o,Q]$.
Notice that $[o,\phi P]=\phi[o,P]$ for every $\phi\in\sln$ and $P\in\MP^n$.
Again by Theorem \ref{thm:contrav}, we obtain that \(Z_2\) is also an $\sln$ contravariant valuation on $\MP^n$.

It remains to show the reverse statement. 
Indeed, by Lemma \ref{lemuq2}, we only need to show that for every $s>0$, $Z$ has the corresponding representation on $sT^d$ for $0\leq d\leq n$ and on $s[e_1,\dots,e_d]$ for $1 \le d \le n$.
By Theorem \ref{thm:tensor:p+2}, there exists a Cauchy function $\eta_1$ such that
\[Z(sT^d)=\sum_{u \in \MN_o (sT^d)} h_{sT^d}^{-p}(u) \eta_1(V(sT^d,u)) \symtp{u}{p}.\]
Let $\spl$ be the set of simplices in $\R^n$ with one vertex at the origin.
For any $T\in\splo$, we write $\tilde T$ as its facet opposite to the origin.
We define the new map $\tilde Z :\spl\to\tenset{n}{p}$ by $\tilde Z(T)= Z (\tilde T)$ for every $T\in\splo$ and $\tilde Z(\{o\})=o$.
It is not hard to check that $\tilde Z $ is an $\sln$ contravariant valuation on $\spl$. 
Since the proofs in \S \ref{sec:n>2} are based only on triangulations in $\spl$  (and not on other polytopes), one can see that Theorem \ref{thm:tensor:p+2} also holds on $\spl$.
Hence there exists a Cauchy function $\eta_2$ such that
\[Z(s\tilde T^d)=\tilde Z(sT^d)=\sum_{u \in \MN_o (sT^d)} h_{sT^d}^{-p}(u) \eta_2(V(sT^d,u)) \symtp{u}{p}.\]
Now, we set $\zeta_1=\eta_1-\eta_2$ and $\zeta_2=\eta_2$. Then
\eqref{eqn:61} holds for $sT^d$ for $0 \le d \le n$ and $s\tilde T^d$ for $1 \le d \le n$ (recall that our cone-volume measure is a signed measure defined in \S \ref{sec:preliminaries}). This completes the proof via Lemma \ref{lemuq2}.
\end{proof}

\section{Tensor valuations of orders larger than $n$.}
In this section, we briefly list all $\sln$ contravariant tensor valuations of order larger than $n$ that we currently know. For brevity, we omit the details of the proof.
\begin{itemize}
\item $\abb{\MPon \ni P \mapsto \sum_{u \in \mathcal {N}_o (P)} h_P^{-p}(u) \zeta(V(P,u)) \symtp{u}{p}} \in \tenval(\MPon;\tenset{n}{p})$, where $\zeta$ is a Cauchy function.

\item $\MPon \ni P \mapsto \int_P x^p \dif x$ is an $\sln$ covariant tensor valuation. It induces a $Z \in \tenval(\MPon;\tenset{n}{(n-1)p})$ defined by
\begin{align}
ZP=\int_P (\cross{x})^p \dif x,
\end{align}
namely, for any $y_1,\dots,y_{(n-1)\times p} \in \R^n$,
\begin{align}
Z P(y_1,\dots,y_{(n-1) \times p})
=\int_P \prod_{j=1}^p \det (y_{(n-1)\times(j-1)+1},\dots,y_{(n-1)\times j},x) \dif x
\end{align}

\item $\MPon \ni P \mapsto \sum_{u \in \mathcal {N}_o (P)}  h_P(u)^{-i+1} \int_{F(P,u)} x^{p-i} \otimes u^i \dif \MH^{n-1} (x)$ for $0 \le i \le p$ is a tensor valuation. It induces a $Z \in \tenval(\MPon;\tenset{n}{(n-1)(p-i)+i})$ defined by
\begin{align}
ZP=\sum_{u \in \mathcal {N}_o (P)}  h_P(u)^{-i+1} \int_{F(P,u)} (\cross{x})^{p-i} \otimes u^i \dif \MH^{n-1} (x).
\end{align}
Remark:  For $i=0$ or $i=p$, these recover the previous two cases. The symmetrization of this tensor was introduced by Haberl and Parapatits \cite{HP17TV} for $P$ containing the origin in its interior.

\item 
$\abb{\MPon \ni P \mapsto (c_0'(-1)^{\dim P}V_0(o \cap \relint P) + c_0V_0(P) + \eta(V_n(P))) \lct^p}$ belongs to \\
$\tenval(\MPon;\tenset{n}{np})$, where $c_0',c_0$ are constants, $\eta$ is a Cauchy function and $\lct^p=\underbrace{\lct \otimes \cdots \otimes \lct}_{p\text{-times}}$.

\item If $Z \in \tenval(\MPon;\tenset{n}{q})$ such that $p-q=rn$. Then
$$P \mapsto Z(P) \otimes \lct^r \in \tenval(\MPon;\tenset{n}{p})$$

\item Only for $n=2$, $\dten{p}(P) \in \tenval(\MPon;\tenset{n}{p})$ is defined in \S \ref{sec:contra}. Let us recall the definition here.
$$\dten{p}(P):=v^p+w^p,$$
if $\dim P=2$ and $P$ has two edges $[o,v]$ and $[o,w]$,
or $\dim P=2$ and $P$ has an edge $[v,w]$ such that $o \in (v,w)$;
$$\dten{p}(P):=2(v^p+w^p),$$
if $\dim P=1$ and $P=[v,w]$ contains the origin; and
$$\dten{p}(P):=0,$$ otherwise.
\end{itemize}

\textbf{Question}: Do these valuations, together with their permutations, span the space of $\sln$ contravariant tensor valuations for $n \ge 2$?

The difficulty of this question can be understood from the perspective of our inductive framework. Our classification strategy proceeds by a double induction on both the dimension and the tensor order, with the planar case and the lower orders serving as the base of the induction. For orders larger than $n$, this planar base case already encounters severe combinatorial and algebraic complexity, and its complete resolution is currently out of reach. The induction step itself reduces the higher--dimensional and higher--order cases to this planar starting point and, in principle, remains valid. The obstruction is therefore not a failure of the induction method, but rather the fact that its base case---the classification in the plane---remains unresolved.

\addcontentsline{toc}{section}{Acknowledgment}
\section*{Acknowledgments}
The author wish to thank the referees for their many valuable remarks that helped to improve the manuscript. The authors are grateful to Thomas Wannerer for his enlightening discussions, especially on the viewpoint of irreducible representations.
The work of the first author is supported in part by the National Natural Science Foundation of China (12571144, 12201388) and
the work of the second author is supported in part by the National Natural Science Foundation of China (12471055)
and the Natural Science Foundation of Gansu Province (23JRRG0001).

\end{document}